\documentclass[a4paper,9pt]{article}
\usepackage[latin1]{inputenc}
\usepackage{graphicx}
\usepackage{amsmath}
\usepackage{amsthm}
\usepackage{amssymb}
\usepackage{color}
\usepackage{mathrsfs}
\usepackage{fancyhdr}
\usepackage{setspace}
\usepackage[center]{subfigure}
\usepackage{float}
\usepackage[center,small,bf,ruled]{caption}
\usepackage{vmargin}
\setmarginsrb{25mm}{25mm}{25mm}{25mm}
             {0mm}{10mm}{0mm}{10mm}

\begin{document}

\numberwithin{equation}{section}

\theoremstyle{thm}
\newtheorem {teo} {Theorem} [section]
\newtheorem {cor} {Corollary} [section]
\newtheorem {lem} {Lemma} [section]

\theoremstyle{definition}
\newtheorem {oss} {Remark} [section]

\newcommand {\Dim} {\textsc{Proof\\ }}
\newcommand {\Fine} { $\blacksquare$ \\}
\date{} 
\title{{ Cascades of Particles Moving at Finite Velocity\\ in Hyperbolic Spaces}}
\maketitle

\long\def\symbolfootnote[#1]#2{\begingroup
\def\thefootnote{\fnsymbol{footnote}}\footnote[#1]{#2}\endgroup} 

\author{
{\center
{ \large V. Cammarota \symbolfootnote[2]{Dipartimento di Statistica, Probabilit\`a e Statistiche applicate, University of Rome `La Sapienza',  P.le Aldo Moro 5, 00185 Rome, Italy.  Tel.: +390649910499, fax: +39064959241.  E-mail address: valentina.cammarota@uniroma1.it.}\;\;\;\;\;\;  E. Orsingher \symbolfootnote[9]{{\it Corresponding author}. Dipartimento di Statistica, Probabilit\`a e Statistiche applicate, University of Rome `La Sapienza',  P.le Aldo Moro 5, 00185 Rome, Italy.  Tel.: +390649910585, fax: +39064959241. E-mail address: enzo.orsingher@uniroma1.it. }} \\
}
}
\vspace{1cm}

\begin{abstract}
A branching process of particles moving at finite velocity over the geodesic lines of the hyperbolic space (Poincar\'e half-plane and Poincar\'e disk) is examined. Each particle can split into two particles only once at Poisson paced times and deviates orthogonally  when splitted. At time $t$, after $N(t)$ Poisson events, there are $N(t)+1$ particles moving along different geodesic lines. We are able to obtain the exact expression of the mean hyperbolic distance of the center of mass of the cloud of particles. We derive such mean hyperbolic distance from two different and independent ways and we study the behavior of the relevant expression as $t$ increases and for different values of the parameters $c$ (hyperbolic velocity of motion) and $\lambda$ (rate of reproduction). The mean hyperbolic distance of each moving particle is also examined and a useful representation, as the distance of a randomly stopped particle moving over the main geodesic line, is presented. 
\end{abstract}

{\small {\bf Keywords}: Branching processes, difference-differential equations, hyperbolic Brownian motion, hyperbolic trigonometry, Laplace transforms, non-Euclidean geometry, random motions.}\\

AMS Classification \\

\section{Introduction}
\hspace{5mm}  Random motions in hyperbolic spaces have been studied since the Fifties and much emphasis has been placed on the so-called hyperbolic Brownian motion on the Poincar\'e  half-plane (see, e.g., Gertsenshtein and Vasiliev \cite{get}, Getoor \cite{getoor}, Gruet \cite{gruet2}, and Lao and Orsingher \cite{Lao}).

 Hyperbolic Brownian motion has been revitalized by mathematical finance since some exotic financial products (Asian options) have a strict connection with the stochastic representation of the hyperbolic Brownian motion  (Yor \cite{MYor}).
 
Branching hyperbolic Brownian motion has been analyzed by Lalley and Sellke \cite{Lalley} who investigated the connection between the birth rate and the underlying dynamics in supercritical and subcritical cases. Also Kelbert and Suhov \cite{kelbert}, \cite{kelbert1} have studied the asymptotic behavior of the hyperbolic branching Brownian motion, developing the ideas in \cite{get} and  \cite{Lalley}.  

 The space on which the above considered hyperbolic Brownian motions develop is the Poincar\'e half-plane (and its higher-dimensional equivalents, see Gruet  \cite{gruet}) or the Klein model (see \cite{kelbert} and \cite{kelbert1}).
 
 The half-plane Poincar\'e model is a fine tool to describe the light propagation in a non-homogeneous medium where, on the basis of Fermat's principle, the angle $\alpha(y)$ between the tangent to the geodesic curve at point of ordinate $y$ satisfies the equality $[{\sin \alpha(y)}]/[{cy}]={1}/{k}$. The randomly scattered irregularities in the medium cause deviations of the trajectories and their decomposition into different rays.   
 
 Random motions with finite velocity have been considered in Orsingher and De Gregorio \cite{OD} (on $H_2^+$ and on the Poincar\'e disk) with the assumption that they develop on geodetic lines and have independent components.
 
 Several models of random motions in $H_2^+$ with finite velocity have been examined in Cammarota and Orsingher \cite{travel}, where the components of the motion have been assumed dependent and the particle moves on mutually orthogonal geodesic lines.
 
Here we study a random motion of a cloud of particles moving at finite velocity on geodesic lines of the hyperbolic space $H_2^+$. Such particles are generated by successive splitting out of a unit-mass particle initially placed at the origin $O$ of $H_2^+$.  The disintegration process of the unit-mass particle is governed by an underlying Poisson process of rate $\lambda$ as follows. The original particle keeps moving on the main geodesic line with constant hyperbolic velocity $c$; at the first Poisson event it breaks into two parts of equal mass $1/2$. One particle continues its motion on the main geodesic line, while the other one deviates orthogonally. At the second Poisson event the deviated particle is separated into two pieces each of mass $1/2^2$; one piece continues its motion on the same geodesic line whereas the other one starts moving on the geodesic line orthogonal to that joining its position with the origin $O$. In general (see Figure \ref{piano}), at the $k$-th Poisson event, only the deviating particle of mass $1/2^{k}$ breaks into two fragments of equal mass $1/2^{k+1}$; the first continues its motion on the same geodesic line while the other one deviates orthogonally.
 
 If up to time $t$, $N(t)$ Poisson events have occurred, we have $N(t)+1$ particles running, at a constant hyperbolic velocity $c$, along different geodesic lines with a mass depending on the instant of separation from the generating particle.
 
  In the above branching process each particle can reproduce only once and the particle splitting at time $T_k$ of the $k$-th Poisson event and which never more disintegrates will preserve its mass, equal to $1/2^k$, for the successive time interval $(T_k, t)$.

 Our main result concerns the dynamics of the center of mass $cm$ of the cloud of particles performing the branching and diffusion process. In particular, we are able to give an exact expression for the mean hyperbolic distance from the origin $O$, $\eta_{cm}(t)$, of the center of mass $cm$ at any time $t>0$
 \begin{eqnarray} \label{12:03}
E\{ \cosh \eta_{cm}(t)\}&=&\frac{2^3 c^2 e^{-\frac{3}{2^2}\lambda t}}{\sqrt{\lambda^2+2^4 c^2}}  \left\{ \frac{e^{-\frac{t}{2^2} \sqrt{\lambda^2+2^4c^2}}}{3 \sqrt{\lambda^2+2^4c^2}+5 \lambda}+\frac{e^{\frac{t}{2^2}\sqrt{\lambda^2+2^4c^2}}}{3 \sqrt{\lambda^2+2^4c^2}-5 \lambda} \right\}  \\
&&+ \frac{\lambda+2c}{2(\lambda+3c)}e^{ct}+\frac{\lambda-2c}{2(\lambda-3c)}e^{-ct}. \nonumber
\end{eqnarray}
We give two different and independent proofs of the above result: our first technique is based on Laplace transforms, while the other one brings about the following non-homogeneous second-order differential equation 
\begin{equation} 
 \frac{\mathrm{d}^2}{\mathrm{d}t^2} u -c^2 u=\frac{\lambda c^2 e^{-\frac{3}{2^2}\lambda t}}{\sqrt{\lambda^2+2^4 c^2}}  \left\{e^{-\frac{t}{2^2} \sqrt{\lambda^2+2^4 c^2}}  -e^{\frac{t}{2^2} \sqrt{\lambda^2+2^4 c^2}}  \right\} \nonumber
\end{equation}
which is satisfied by (\ref{12:03}). 

The behavior of the hyperbolic distance of each individual particle can be compared with result (\ref{12:03}). In a previous paper (see \cite{travel}) we have shown that the mean hyperbolic distance $\eta(t)$ of the particle which underwent changes of direction at all Poisson events is   
\begin{eqnarray} \label{aitia}
E\{ \cosh \eta(t)\} 
&=&\frac{2 c^2 e^{-\frac{\lambda t}{2}}}{\sqrt{\lambda^2+2^2c^2}}  \left\{ \frac{e^{-\frac{t}{2}\sqrt{\lambda^2+2^2c^2}}}{\sqrt{\lambda^2+2^2c^2}+\lambda} +  \frac{e^{\frac{t}{2}\sqrt{\lambda^2+2^2c^2}}}{\sqrt{\lambda^2+2^2c^2}-\lambda}  \right\}. 
\end{eqnarray}
If in (\ref{aitia}) the Poisson rate $\lambda$ is replaced by $\lambda/2$ we see that the exponential terms in (\ref{12:03}) and (\ref{aitia}) have the same form but different weights. 

We also examine the mean hyperbolic distance of each individual particle which stops changing direction after the $k$-th Poisson event. Our main result shows that 
\begin{equation} \label{18:13}
E\{\cosh \eta_k(t) I_{\{N(t)\ge k   \}}\}= \frac{1}{2^k} \int_{0}^{t} \cosh c(t-s) h(k,c,s) g(s;k,\lambda) \; \mathrm{d} s
\end{equation}
 where $g(s;k,\lambda)=\frac{e^{-\lambda s} \lambda^k s^{k-1}}{\Gamma(k)}$ is a Gamma distribution and $h(k,c,s)= \sum_{r=0}^{k} \binom{k}{r}  E_{Y_{r,k}} \{e^{cs(2Y_{r,k}-1)} \}$ with $Y_{r,k} \sim \mathrm{Beta}(r,k-r)$ so that $(2Y_{r,k}-1)\in(-1,1)$ and with the assumption that $Y_{0,k}=1$ and $Y_{k,k}=-1$. Thus (\ref{18:13}) shows that the mean hyperbolic distance, at time $t$, of a particle generated at the $k$-th Poisson event can be seen as the mean hyperbolic distance of a particle which never deviates from the main geodesic line and which starts moving at a random time with law  $h(k,c,s) g(s;k,\lambda)$.

 \section{Some geometrical features of the hyperbolic spaces}

We present in this section some basic features of the Poincar\'e half-plane $H_2^+=\{ (x,y):\, y>0\}$ which is endowed with the metric
\begin{equation} \label{metric}
\mathrm{d}s=\frac{\sqrt{(\mathrm{d}x)^2+(\mathrm{d}y)^2}}{y}.
\end{equation}
Some informations on hyperbolic spaces and non-Euclidean geometry can be found in Faber \cite{Faber} and Meschkowski \cite{wski}. 
The position of points in $H_2^+$ can be given in Cartesian coordinates $(x, y)$ or in hyperbolic coordinates $(\eta, \alpha)$. These are connected by means of the well-known relationships 
\begin{equation}\label{rel}
\left\{
\begin{array}{lr} x=\frac{\sinh \eta \cos \alpha}{\cosh \eta - \sinh \eta \sin \alpha}, &  \eta>0,\\
y=\frac{1}{\cosh \eta - \sinh \eta \sin \alpha}, &  -\frac{\pi}{2}<\alpha<\frac{\pi}{2}.
\end{array}
\right.
\end{equation}
The hyperbolic coordinate $\eta$ represents the distance of $(x, y)$ from the origin $O=(0,1)$ of $H_2^+$ measured by means of  the metric (\ref{metric}). The coordinate $\alpha$ is the angle of the tangent in $O$ to the half-circumference joining $O$ with $(x,y)$ (see Rogers and Williams \cite{rog} page 213 and Cammarota and Orsingher \cite{travel} for some details). \\
For each point $(x,y) \in H_2^+$ it is possible to obtain the hyperbolic distance $\eta$ as well as the angle $\alpha$ by means of the formulas 
\begin{equation} \label{rell}
\cosh \eta=\frac{x^2+y^2+1}{2y}, \hspace{1cm} \tan \alpha=\frac{x^2+y^2-1}{2x}.
\end{equation}
Formulas (\ref{rell}) are easily derived from (\ref{rel}) (see  \cite{OD}). We can obtain formulas (\ref{rel}) from (\ref{rell}) as follows. By substituting 
\begin{equation}
x=\tan \alpha \pm \sqrt{\tan ^2 \alpha +1-y^2}, \nonumber
\end{equation}
in the first relationship of (\ref{rell}) we have that
\begin{equation}
y \cosh \eta -1- \tan ^2 \alpha =\pm \tan \alpha \sqrt{\tan ^2 \alpha+1-y^2}, \nonumber
\end{equation}
and after some manipulations we arrive at 
\begin{eqnarray}
0= (y \cosh \eta-1)^2-y^2 \sin ^2 \alpha \sinh^2 \eta= ( y \cosh \eta-1-y \sin \alpha \sinh \eta)( y \cosh \eta-1+y \sin \alpha \sinh \eta), \nonumber
\end{eqnarray}
which yields the second formula of (\ref{rel}). Since from  (\ref{rell}) $y \cosh \eta=x \tan \alpha-1,$ also the second relationship of (\ref{rel}) immediately follows.    

In $H_2^+$ the trigonometrical formulas we need are either the Pythagorean theorem for right triangles 
\begin{equation} \label{pit}
\cosh \eta=\cosh \eta_1 \cosh \eta_2,
\end{equation}
or its Carnot extension for arbitrary triangles
\begin{equation} \label{carnot}
\cosh \eta=\cosh \eta_1 \cosh \eta_2 -\sinh \eta_1 \sinh \eta_2 \cos(\alpha_1-\alpha_2).
\end{equation}
Clearly, if $\alpha_1-\alpha_2=\pi/2$, formula (\ref{carnot}) becomes (\ref{pit}).

The half-plane $H_2^+$ can be mapped onto the disk $D=\{(u, v): u^2+v^2<1\}$ by means of the conformal mapping 
\begin{equation} \label{11:47}
w=\frac{iz+1}{z+i}.
\end{equation} 
The $x$-axis of $H_2^+$ is mapped onto the boundary $\partial D$ of $D$ while the origin $O$ is transformed into the center of $D$. An arbitrary point $(x, y) \in H_2^+$ is mapped onto the point $(u, v) \in D$ with coordinates 
\begin{equation}
u=\frac{2x}{x^2+(y+1)^2}, \hspace{2cm} v=\frac{x^2+y^2-1}{x^2+(y+1)^2}. \nonumber
\end{equation} 
A point $(u,v) \in D$ is instead mapped by (\ref{11:47}) into the point $(x, y)$ with coordinates 
\begin{equation}
x=\frac{2 u}{u^2+(1-v)^2}, \hspace{2cm} y=\frac{1-(u^2+v^2)}{u^2+(1-v)^2}. \nonumber
\end{equation}
A similar mapping is the so-called Cayley transformation which reads 
\begin{equation} \label{11:48}
w=\frac{i-z}{i+z},
\end{equation}
and it slightly differs from (\ref{11:47}). The geodesic lines of $H_2^+$ with radius $r$ and center at $(x_0,0)$, are mapped by (\ref{11:47}) into arcs of circumferences inside $D$ with center at 
\begin{equation}\left( \frac{2 x_0}{x_0^2-r^2-1}, \;\; \frac{x_0^2-r^2-1}{x_0^2-r^2+1} \right) \nonumber \end{equation}
and with radius $R$ given by $R^2=\frac{4 r^2}{(x_0^2-r^2-1)^2}$.

\section{Description of the randomly moving and branching model}

\begin{figure}[h]
 \centering
    {\includegraphics[width=15.5cm, height=10.7cm]{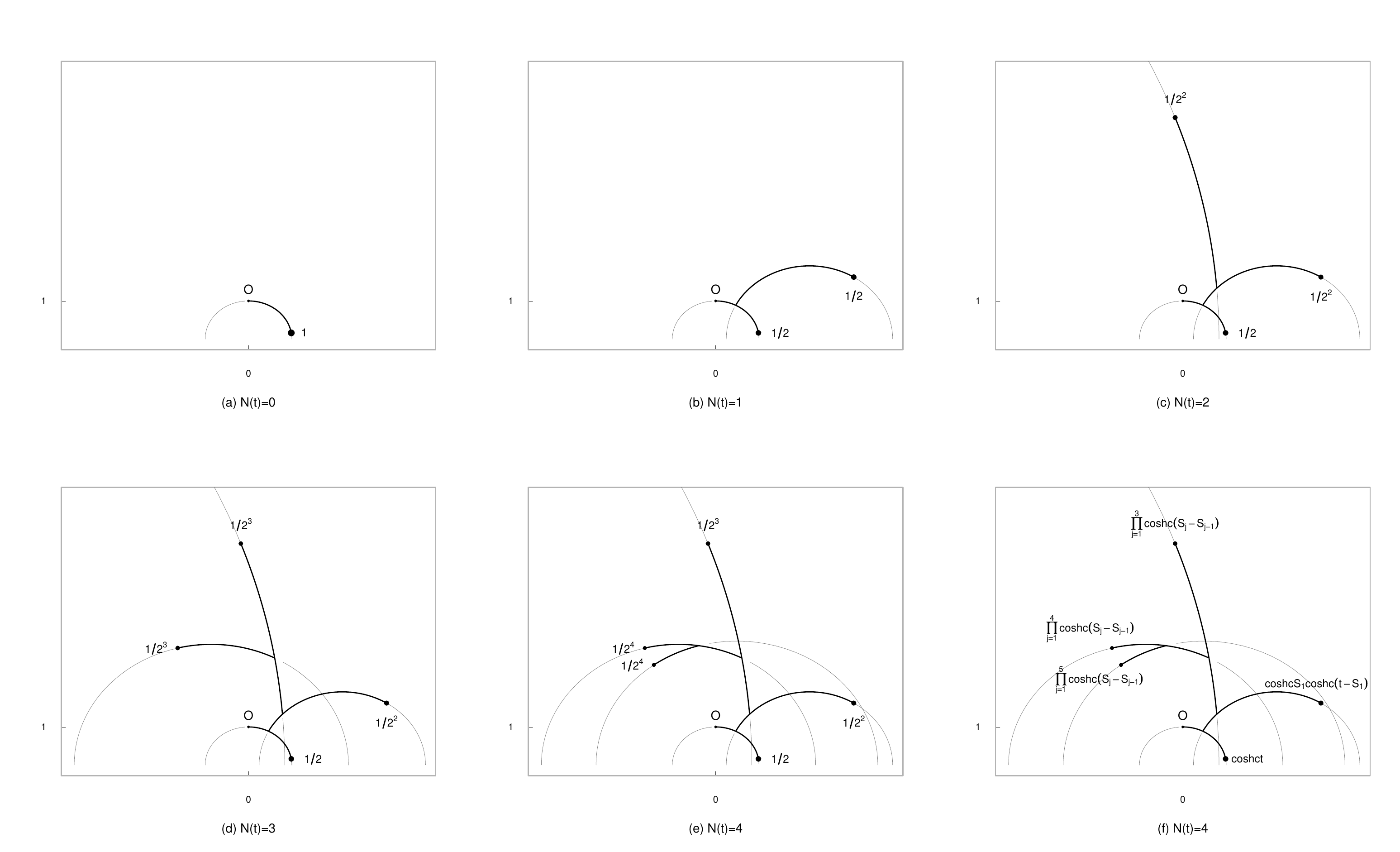}}
 \caption{In (a), the trajectory of the unit-mass particle initially placed at the origin $O$ of $H_2^+$ and moving on the main geodesic line is shown. When $N(t)=0$, no disintegration occurs. In (b), (c), (d), and (e), the trajectories of the particles generated by $N(t)=1,2,3,4$ Poisson events are plotted. The relevant mass associated with each particle is also indicated by a suitable label. In (f), for each particle the relevant hyperbolic cosine of the hyperbolic distance from the origin is also reported.} \label{piano}
 \end{figure}

We assume that a unit-mass particle is placed at time $t=0$ at the origin $O$ of $H_2^+$ and starts moving on the main geodesic line represented in $H_2^+$ by the half-circle of radius $1$ passing through $O$. This particle chooses with probability $1/2$ one of the two possible directions and moves with constant hyperbolic velocity equal to $c$ (see Figure \ref{piano} (a)). The hyperbolic velocity
\begin{equation}
c=\frac{\mathrm{d}s}{\mathrm{d}t}=\frac{1}{y}\sqrt{\left(\frac{\mathrm{d}x}{\mathrm{d}t}\right)^2+\left(\frac{\mathrm{d}y}{\mathrm{d}t}\right)^2} \nonumber
\end{equation}
is assumed to be constant. For an Euclidean observer, the closer to the $x$-axis is the moving particle the slower it moves.

A Poisson process of rate $\lambda$ governs the changes of direction. At the first Poisson event the particle splits into two pieces of equal mass: one continues its motion on the same geodesic line while the other one starts moving (in one of the two possible directions) on the geodesic line orthogonal to the previous one (see Figure \ref{piano} (b)). In general, at the $k$-th Poisson event, the deviating particle of mass $1/2^k$ undergoes a further decomposition: one splinter of mass $1/2^{k+1}$ continues undisturbed its motion, while the other one, also of mass $1/2^{k+1}$, is forced to move onto the geodetic line orthogonal to that joining $O$ with the position it occupied at the time where the splitting took place. 

 Therefore, if no Poisson event occurs (i.e., $\{N(t)=0\}$, where $N(t)$ is the number of Poisson events in $[0,t]$), the unit-mass particle is located, at time $t$, on the first geodesic line at an hyperbolic distance from $O$ equal to $\eta_0(t)=ct$ (see Figure \ref{piano} (a)).  

If one Poisson event  happens, at time $S_1<t$ (i.e., $\{N(t)=1\}$), then, at time $t$, one fragment of mass $1/2$ will be at distance $\eta_0(t)=ct$ on the first geodesic line, while the other splinter, also of mass $1/2$, will be located at hyperbolic distance $\eta_1(t)$ given by
\begin{equation} \label{16:44}
\cosh \eta_1(t)=\cosh c\,S_1 \cosh c(t-S_1),
\end{equation}  
where in formula (\ref{16:44}) the hyperbolic Pythagorean theorem (\ref{pit}) has been applied (see Figure \ref{piano} (b)). 

If $N(t)=n$, the process described above produces $n+1$ splinters. The particle generated at the $k$-th Poisson event, $k=0, \cdots n-1$, (and which will not break up after the $k$-th event) has mass $1/2^{k+1}$ and is located, at time $t$, at the hyperbolic distance $\eta_k(t)$ from $O$. Such a distance, in force of the hyperbolic Pythagorean theorem, reads  
\begin{equation}
\cosh \eta_k(t)=\prod_{j=1}^{k+1} \cosh c(S_j-S_{j-1}),    \nonumber            
\end{equation}
where $S_0=0$, $S_{k+1}=t$ and the random times $S_j$ with $j=1, \cdots k$ represent the instants where the Poisson events happen and the deviations of motion occur. The last splinter, which has changed direction at all fission events, has mass $1/2^n$ and, at time $t$, is at an hyperbolic distance from $O$ given by
\begin{equation}
\cosh \eta_n (t)=\prod_{j=1}^{n+1} \cosh c(S_j-S_{j-1}). \nonumber
\end{equation}
where $S_0=0$ and $S_{n+1}=t$ (see Figure \ref{piano} (f)). 

In general, if the number of splits recorded is $N(t)$, the number of particles is $N(t)+1$ and each runs on a different  geodesic line of $H_2^+$. The hyperbolic distance of the center of mass $cm$ of the cloud of moving particles, at time $t>0$, is denoted by $\eta_{cm}(t)$ and is represented by 
\begin{equation} \label{dcm}
\cosh \eta_{cm}(t)=\sum_{k=0}^{N(t)-1} \frac{1}{2^{k+1}} \prod_{j=1}^{k+1} \cosh c(S_j-S_{j-1}) I_{\{N(t)>0\}}+ \frac{1}{2^{N(t)}} \prod_{j=1}^{N(t)+1} \cosh c(S_j-S_{j-1}),
\end{equation} 
where $S_1, S_2\dots S_{N(t)}$ are the random times at which the disintegrations occur, and $S_0=0$, $S_{N(t)+1}=t$. The second term in (\ref{dcm}) refers to the splinter which underwent all disintegrations occurred until time $t$, while the first one is related to those particles produced during the branching process. 

The assumption that the particles deviate on geodesic lines orthogonal to those joining the origin $O$ with their current position is crucial since it makes the hyperbolic Pythagorean theorem applicable. Otherwise it should be necessary to apply the Carnot hyperbolic formula and this would make the analytic treatment of the problem extremely difficult. 

Under the condition that $N(t)=n$, the mean hyperbolic distance of the center of mass of the cloud of  particles at time $t$ is 
 \begin{eqnarray} \label{11:17}
 E \{\cosh \eta_{cm}(t) | N(t)=n\}&=&\sum_{k=0}^{n-1} \frac{1}{2^{k+1}} \prod_{j=1}^{k+1} \cosh c(S_j-S_{j-1}) I_{\{n>0\}} +\frac{1}{2^n}\prod_{j=1}^{n+1} \cosh c(S_j-S_{j-1}).\;\;\;
 \end{eqnarray} 
We observe that the $n$ instants $S_1, \cdots, S_n$ where the fissions take place are uniformly distributed under the condition that $N(t)=n$, and possess density 
\begin{equation}
Pr\{S_1 \in \mathrm{d} s_1,\cdots, S_n \in \mathrm{d}s_n\}= \frac{n!}{t^n} \mathrm{d}s_1 \cdots \mathrm{d}s_n \nonumber
\end{equation}
for $0=s_0<s_1< \cdots <s_{n+1}=t$. Therefore, under the condition that the number of splitting events is $N(t)=n$, the hyperbolic distance $\eta_k(t)$ of the $k$-th splinter is for $k=0, \cdots n-1$
\begin{eqnarray} \label{Gk}
E\{ \cosh \eta_k(t)|N(t)=n \}&=& \frac{n!}{t^n} \int_{0}^{t} \mathrm{d}s_1 \cdots \int_{s_{k-1}}^{t} \mathrm{d}s_k \cdots \int_{s_{n-1}}^{t} \mathrm{d}s_n \prod_{j=1}^{k+1} \cosh c(s_j -s_{j-1})\nonumber \\ &=& \frac{n!}{t^n} \int_{0}^{t} \mathrm{d}s_1 \cdots \int_{s_{k-1}}^{t} \mathrm{d}s_k \frac{(t-s_k)^{n-k}}{(n-k)!} \prod_{j=1}^{k+1} \cosh c(s_j-s_{j-1})  \nonumber \\ &=& \frac{n!}{t^n}  G_{n,k}(t),
\end{eqnarray}
and for $k=n$
\begin{eqnarray} \label{Gn}
E\{\cosh \eta_n(t)|N(t)=n\}&=& \frac{n!}{t^n} \int_{0}^{t} \mathrm{d}s_1 \cdots \int_{s_{n-1}}^{t} \mathrm{d}s_n  \prod_{j=1}^{n+1} \cosh c(s_j -s_{j-1}) \nonumber \\
 &=& \frac{n!}{t^n}  G_{n,n}(t).
\end{eqnarray}

\begin{figure}[h] \label{Disco}
 \centering
   {\includegraphics[width=15.5cm, height=11.2cm]{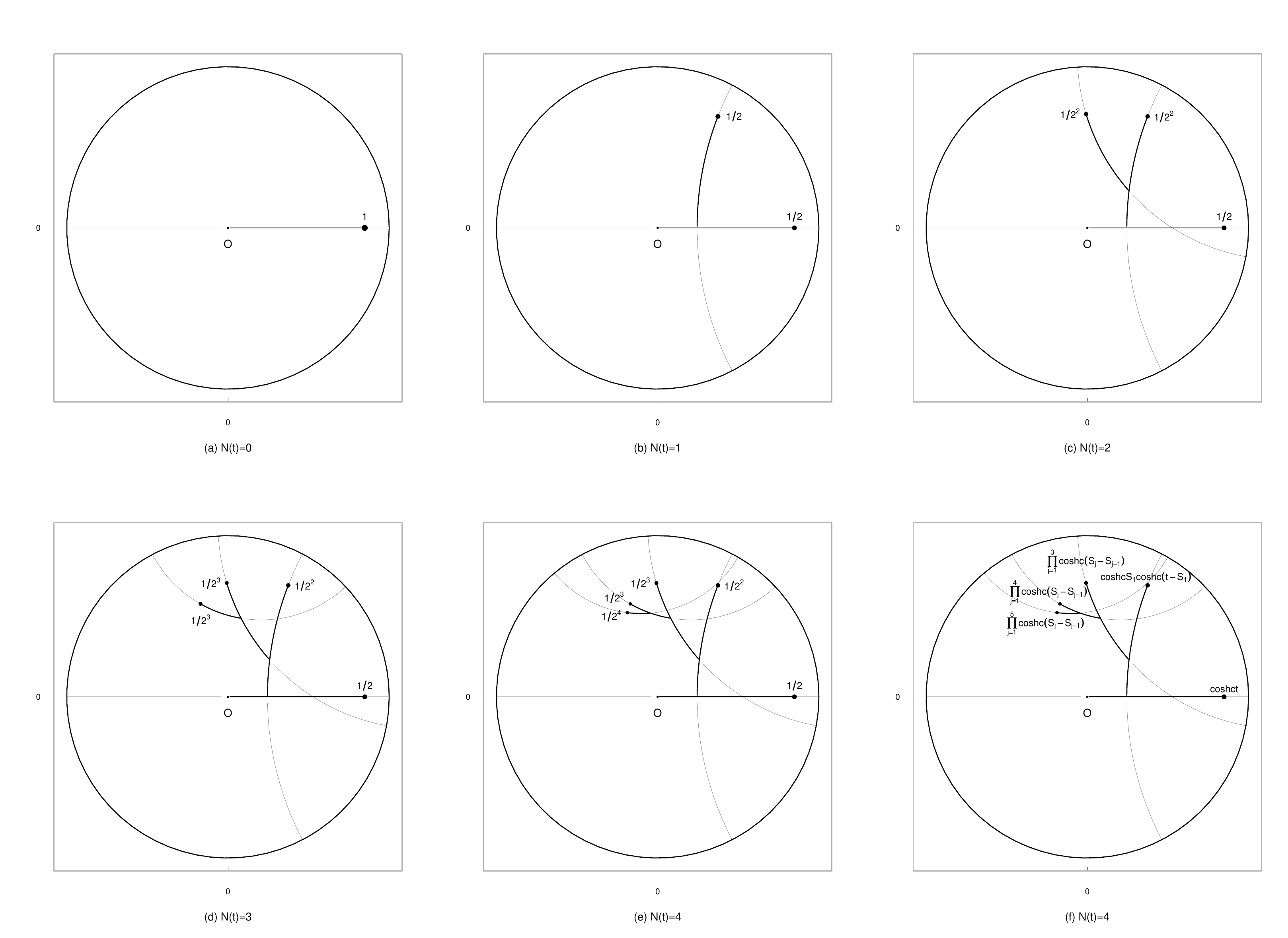}}
 \caption{Same trajectories as in Figure \ref{piano} represented through the Poincar\'e disk model.} \label{disco}
 \end{figure}

The branching process described above can be adapted to the Poincar\'e disk $D$: since the mappings (\ref{11:47})  and (\ref{11:48}) preserve the hyperbolic distance, the trajectories of splitting and moving particles can be conveniently depicted in $D$ as in $H_2^+$ (see Figure \ref{piano} and \ref{disco}).

We restrict ourselves to the mean hyperbolic distance because this leads to fine explicit results. The analysis of the distribution of the hyperbolic distance, even in the case of a random motion of an individual particle that changes direction at all Poisson events, implies a much more complicated analysis and is almost intractable since multiple integrals of the form 
\begin{equation}
E\{e^{i \alpha \cosh \eta (t)}| N(t)=n\}=\frac{n!}{t^n} \int_{0}^{t} \mathrm{d}s_1 \cdots \int_{s_{n-1}}^{t}   e^{i \alpha \prod_{j=1}^{n+1} \cosh c(s_j-s_{j-1})} \mathrm{d}s_n \nonumber
\end{equation}      
must be evaluated (see Cammarota and Orsingher \cite{travel} for details on this point).  

\section{Mean hyperbolic distance of the system of randomly moving and disintegrating particles: the Laplace transform approach}

We are able to obtain the explicit form of the mean-value of (\ref{dcm}) by means of two different and independent approaches. In this section we present the Laplace-transform derivation which leads to our first theorem.   

\begin{teo}
The mean-value of the hyperbolic distance (\ref{dcm}) of the center of mass is 
\begin{eqnarray} \label{12:36}
E\{ \cosh \eta_{cm}(t)\}&=&\frac{2^3 c^2 e^{-\frac{3}{2^2}\lambda t}}{\sqrt{\lambda^2+2^4 c^2}}  \left\{ \frac{e^{-\frac{t}{2^2} \sqrt{\lambda^2+2^4c^2}}}{3 \sqrt{\lambda^2+2^4c^2}+5 \lambda}+\frac{e^{\frac{t}{2^2}\sqrt{\lambda^2+2^4c^2}}}{3 \sqrt{\lambda^2+2^4c^2}-5 \lambda} \right\}  \\
&&+ \frac{\lambda+2c}{2(\lambda+3c)}e^{ct}+\frac{\lambda-2c}{2(\lambda-3c)}e^{-ct}, \hspace{3cm}t>0. \nonumber
\end{eqnarray}
\end{teo}  
\Dim
In view of (\ref{11:17}), (\ref{Gk}) and (\ref{Gn}), we first note that 
\begin{eqnarray}
E\{\cosh \eta_{cm}(t)| N(t)=n\}&=&\frac{n!}{t^n} \sum_{k=0}^{n-1} \frac{1}{2^{k+1}} \int_{0}^{t} \mathrm{d}s_1 \cdots \int_{s_{n-1}}^{t} \mathrm{d}s_{n} \prod_{j=1}^{k+1} \cosh c(s_j-s_{j-1}) \nonumber \\
&&+\frac{n!}{t^n} \frac{1}{2^n} \int_{0}^{t} \mathrm{d}s_1 \cdots \int_{s_{n-1}}^{t} \mathrm{d}s_n \prod_{j=1}^{n+1} \cosh c(s_j-s_{j-1}) \nonumber \\
&=&\left\{
\begin{array}{lr}  \frac{n!}{t^n} \sum_{k=0}^{n-1} \frac{1}{2^{k+1}} G_{n,k}(t)+ \frac{n!}{t^n} \frac{1}{2^n} G_{n,n}(t), \hspace{1cm} n \ge1, \nonumber \\
\\
G_{0,0}(t), \hspace{5.5cm}  n=0.
\end{array}
\right.
\end{eqnarray}
Our task is therefore to study the Laplace transform, 
\begin{eqnarray} \label{co}
\int_{0}^{\infty} e^{-\mu t} E\{\cosh \eta_{cm}(t)\} \mathrm{d}t&=&\int_{0}^{\infty} e^{- \mu t} \sum_{n=0}^{\infty} E\{\cosh \eta_{cm}(t)| N(t)=n\} Pr\{N(t)=n\} \mathrm{d}t \nonumber \\
&=&\sum_{n=1}^{\infty} \lambda^n \sum_{k=0}^{n-1} \frac{1}{2^{k+1}} \int_{0}^{\infty} e^{-(\lambda+\mu )t} G_{n,k}(t) \mathrm{d}t \nonumber \\
&&+ \sum_{n=0}^{\infty} \frac{\lambda^n}{2^n} \int_{0}^{\infty} e^{-(\lambda+\mu)t} G_{n,n} (t) \mathrm{d}t,
\end{eqnarray}
where $\mu>0$. We evaluate the Laplace transform appearing in (\ref{co}) in the following way. If $\gamma=\lambda+\mu$ and $c<\gamma$, then we have, by successively inverting the inner integrals, that 
\begin{eqnarray}
\lefteqn{\int_{0}^{\infty} e^{-\gamma t}G_{n,k}(t) \mathrm{d}t} \nonumber \\
&=&\int_{0}^{\infty} e^{-\gamma t} \left\{\ \int_{0}^{t} \mathrm{d}s_1 \cdots \int_{s_{k-1}}^{t} \mathrm{d}s_k \frac{(t-s_k)}{(n-k)!}^{n-k}\;\; \prod_{j=1}^{k+1} \cosh c(s_j-s_{j-1}) \right\} \mathrm{d}t \nonumber \\
&=&\int_{0}^{\infty} \mathrm{d}s_1 \int_{s_1}^{\infty} e^{-\gamma t} \mathrm{d}t \left\{ \int_{s_1}^{t} \mathrm{d}s_2 \cdots \int_{s_{k-1}}^{t} \mathrm{d}s_k  \frac{(t-s_k)}{(n-k)!}^{n-k} \;\; \prod_{j=1}^{k+1} \cosh c(s_j-s_{j-1}) \right\} \nonumber \\
&=& \int_{0}^{\infty} \mathrm{d}s_1 \int_{s_1}^{\infty} \mathrm{d}s_2 \cdots \int_{s_{k-1}}^{\infty} \mathrm{d}s_k \prod_{j=1}^{k} \cosh c(s_j-s_{j-1})\left\{ \int_{s_k}^{\infty} e^{-\gamma t}  \frac{(t-s_k)}{(n-k)!}^{n-k} \cosh c(t-s_k) \mathrm{d t}\right\} \nonumber \\
&=& \int_{0}^{\infty} \mathrm{d}s_1  \cdots \int_{s_{k-2}}^{\infty} \mathrm{d}s_{k-1} \prod_{j=1}^{k-1} \cosh c(s_j-s_{j-1}) \int_{s_{k-1}}^{\infty} e^{-\gamma s_k} \cosh c(s_k-s_{k-1}) \mathrm{d}s_k  \nonumber \\
&&\times \int_{0}^{\infty} e^{-  \gamma w} \frac{w^{n-k}}{(n-k)!} \cosh c w \; \mathrm{d}w \nonumber \\
&=& \left( \int_{0}^{\infty} e^{-\gamma w} \cosh cw \right)^k \int_{0}^{\infty} e^{-\gamma w} \frac{w^{n-k}}{(n-k)!} \cosh c w\; \mathrm{d}w. \nonumber
\end{eqnarray}
In the last step above the change of variable $s_j-s_{j-1}=w_j$ applied $k$ times leads to the final expression. For $k=n$, from the previous calculations, we obtain 
\begin{equation}
\int_{0}^{\infty} e^{- \gamma t} G_{n,n}(t) \mathrm{d}t= \left( \int_{0}^{\infty} e^{- \gamma w} \cosh cw \;\mathrm{d}w \right)^{n+1}. \nonumber
\end{equation}
Since 
\begin{eqnarray}
\int_{0}^{\infty} e^{-\gamma w} \cosh cw \; \mathrm{d}w&=&\frac{\gamma}{\gamma^2-c^2}, \nonumber \\
\int_{0}^{\infty} e^{-\gamma w} \frac{w^{n-k}}{(n-k)!} \cosh cw \; \mathrm{d}w&=&\frac{1}{2} \left[ \frac{1}{(\gamma-c)^{n-k+1}}+\frac{1}{(\gamma+c)^{n-k+1}} \right], \nonumber
\end{eqnarray}
we have that 
\begin{eqnarray} \label{18:18cat}
\lefteqn{\int_{0}^{\infty} e^{-\mu t} E\{\cosh \eta_{cm}(t)\} \mathrm{d}t} \nonumber \\
&=&\sum_{n=1}^{\infty} \lambda^n \sum_{k=0}^{n-1} \frac{1}{2^{k+1}} \left( \frac{\gamma}{\gamma^2-c^2} \right)^k \frac{1}{2} \left[ \frac{1}{(\gamma-c)^{n-k+1}}+\frac{1}{(\gamma+c)^{n-k+1}} \right] \nonumber \\
&&+ \sum_{n=0}^{\infty} \left(  \frac{\lambda}{2} \right)^n \left( \frac{\gamma}{\gamma^2-c^2} \right)^{n+1} \nonumber \\
&=& \frac{1}{2^2} \sum_{n=1}^{\infty} \lambda^n \left\{ \frac{1}{(\gamma-c)^{n+1}} \sum_{k=0}^{n-1}  \left[\frac{\gamma}{2(\gamma+c)} \right]^k +\frac{1}{(\gamma+c)^{n+1}} \sum_{k=0}^{n-1} \left[\frac{\gamma}{2(\gamma-c)} \right]^k  \right\} \nonumber \\
&&+ \frac{2 \gamma}{2 \gamma^2-2c^2-\gamma \lambda},
\end{eqnarray}
where the last sum converges if $\mu$ satisfies the inequality $2 c^2< \lambda^2+2 \mu^2+3 \lambda \mu$. The double sum in (\ref{18:18cat}) can be calculated by inverting the order of summation in the following way:
\begin{eqnarray} \label{14:12}
\lefteqn{\sum_{n=1}^{\infty} \lambda^n \left\{ \frac{1}{(\gamma-c)^{n+1}} \sum_{k=0}^{n-1}  \left[\frac{\gamma}{2(\gamma+c)} \right]^k +\frac{1}{(\gamma+c)^{n+1}} \sum_{k=0}^{n-1} \left[\frac{\gamma}{2(\gamma-c)} \right]^k  \right\}} \nonumber \\
&=&\frac{1}{\gamma-c} \sum_{k=0}^{\infty} \left[ \frac{\gamma}{2(\gamma+c)} \right]^k \sum_{n=k+1}^{\infty} \left( \frac{\lambda}{\gamma-c} \right)^n+ \frac{1}{\gamma+c} \sum_{k=0}^{\infty} \left[ \frac{\gamma}{2(\gamma-c)} \right]^k \sum_{n=k+1}^{\infty} \left( \frac{\lambda}{\gamma+c} \right)^n   \nonumber \\
&=&\frac{1}{\gamma-c} \sum_{k=0}^{\infty} \left[ \frac{\gamma}{2(\gamma+c)} \right]^k  \left( \frac{\lambda}{\gamma-c} \right)^{k+1}  \sum_{r=0}^{\infty} \left( \frac{\lambda}{\gamma-c} \right)^r \nonumber \\
&&+ \frac{1}{\gamma+c} \sum_{k=0}^{\infty} \left[ \frac{\gamma}{2(\gamma-c)} \right]^k  \left( \frac{\lambda}{\gamma+c} \right)^{k+1}  \sum_{r=0}^{\infty} \left( \frac{\lambda}{\gamma+c} \right)^r \nonumber \\
&=& \frac{\lambda}{(\gamma-c)^2} \sum_{k=0}^{\infty} \left[ \frac{\gamma \lambda }{2(\gamma^2-c^2)} \right]^k \frac{\gamma-c}{\gamma-c-\lambda}  + \frac{\lambda}{(\gamma+c)^2} \sum_{k=0}^{\infty} \left[ \frac{\gamma \lambda }{2(\gamma^2-c^2)} \right]^k \frac{\gamma+c}{\gamma+c-\lambda} \nonumber \\
&=&\lambda  \left[\frac{1}{\gamma-c}\; \frac{1}{\gamma-c-\lambda} + \frac{1}{\gamma+c}\; \frac{1}{\gamma+c-\lambda} \right] \sum_{k=0}^{\infty} \left[ \frac{\gamma \lambda }{2(\gamma^2-c^2)} \right]^k \nonumber \\
&=& \lambda  \left[\frac{1}{\gamma-c}\; \frac{1}{\gamma-c-\lambda} + \frac{1}{\gamma+c}\; \frac{1}{\gamma+c-\lambda} \right] \frac{2(\gamma^2-c^2)}{2 \gamma^2-2c^2- \gamma \lambda}.
\end{eqnarray} 
The inversion of the Laplace transform is made possible by suitably rearranging the expression (\ref{14:12}) as follows: 
\begin{eqnarray} \label{18:30cat}
\lefteqn{\int_{0}^{\infty} e^{-\mu t} E\{\cosh \eta_{cm}(t)\} \mathrm{d}t }\nonumber \\
&=&\frac{\lambda}{2} \left[\frac{1}{\gamma-c}\; \frac{1}{\gamma-c-\lambda}  + \frac{1}{\gamma+c}\; \frac{1}{\gamma+c-\lambda} \right] \frac{\gamma^2-c^2}{2 \gamma^2-2c^2-\gamma \lambda} +\frac{2 \gamma}{ 2 \gamma^2-2c^2-\gamma \lambda} \nonumber\\
&=&\frac{\lambda}{2} \left[ \frac{1}{\lambda+\mu-c}\; \frac{1}{\mu-c}+ \frac{1}{\lambda+\mu+c}\; \frac{1}{\mu+c}  \right] \frac{(\lambda+\mu)^2 -c^2}{\lambda^2+2 \mu^2+3 \lambda \mu -2c^2} +\frac{2\lambda+2\mu}{\lambda^2+2\mu^2+3\lambda \mu-2c^2} \nonumber\\
&=&\frac{\lambda}{2} \left[ \frac{\lambda+\mu+c}{\mu-c}+ \frac{\lambda+\mu-c}{\mu+c}  \right] \frac{1}{\lambda^2+2 \mu^2+3 \lambda \mu -2c^2}+\frac{2\lambda+2\mu}{\lambda^2+2\mu^2+3\lambda \mu-2c^2}.
\end{eqnarray}
By means of the decomposition
\begin{equation}
2\mu^2+3\lambda \mu+\lambda^2-2c^2=2\left[ \mu+\frac{3\lambda-\sqrt{\lambda^2+2^4c^2}}{2^2} \right] \left[ \mu+\frac{3\lambda+\sqrt{\lambda^2+2^4c^2}}{2^2} \right], \nonumber \end{equation}
the expression (\ref{18:30cat}) can be further worked out by writing 
\begin{eqnarray}  \label{18:55}
\lefteqn{\int_{0}^{\infty} e^{-\mu t} E\{\cosh \eta_{cm}(t)\} \mathrm{d}t} \nonumber \\
&=&\frac{\lambda}{2^2} \frac{1}{\left[ \mu+\frac{3\lambda-\sqrt{\lambda^2+2^4c^2}}{2^2} \right] \left[ \mu+\frac{3\lambda+\sqrt{\lambda^2+2^4c^2}}{2^2} \right] }\left[  2 +\frac{\lambda+2c}{\mu-c} +\frac{\lambda-2c}{\mu+c}  \right] \nonumber \\
&&+\frac{\frac{3}{2}\lambda+\frac{\lambda}{2}+2\mu}{2\left[ \mu+\frac{3\lambda-\sqrt{\lambda^2+2^4c^2}}{2^2} \right] \left[ \mu+\frac{3\lambda+\sqrt{\lambda^2+2^4c^2}}{2^2} \right]} \nonumber \\
&=& \left[ \frac{1}{ \mu+\frac{3\lambda-\sqrt{\lambda^2+2^4c^2}}{2^2}}-\frac{1}{ \mu+\frac{3\lambda+\sqrt{\lambda^2+2^4c^2}}{2^2}} \right] \frac{\lambda}{2^2} \frac{2}{\sqrt{\lambda^2+2^4c^2}}  \left[  2 +\frac{\lambda+2c}{\mu-c} +\frac{\lambda-2c}{\mu+c}  \right] \nonumber \\
&&+  \left[ \frac{1}{ \mu+\frac{3\lambda-\sqrt{\lambda^2+2^4c^2}}{2^2}}+\frac{1}{ \mu+\frac{3\lambda+\sqrt{\lambda^2+2^4c^2}}{2^2}} \right]     \frac{\frac{3}{2}\lambda+\frac{\lambda}{2}+2\mu}{2\left( 2\mu+\frac{3}{2}\lambda \right) } \nonumber \\
&=&  \left[ \frac{1}{ \mu+\frac{3\lambda-\sqrt{\lambda^2+2^4c^2}}{2^2}}-\frac{1}{ \mu+\frac{3\lambda+\sqrt{\lambda^2+2^4c^2}}{2^2}} \right] \frac{\lambda}{\sqrt{\lambda^2+2^4c^2}} +  \frac{1}{2} \left[ \frac{1}{ \mu+\frac{3\lambda-\sqrt{\lambda^2+2^4c^2}}{2^2}}+\frac{1}{ \mu+\frac{3\lambda+\sqrt{\lambda^2+2^4c^2}}{2^2}} \right]\nonumber \\
&&+  \left[ \frac{1}{ \mu+\frac{3\lambda-\sqrt{\lambda^2+2^4c^2}}{2^2}}-\frac{1}{ \mu+\frac{3\lambda+\sqrt{\lambda^2+2^4c^2}}{2^2}} \right] \frac{\lambda}{2\sqrt{\lambda^2+2^4c^2}} \left[ \frac{\lambda+2c}{\mu-c}+\frac{\lambda-2c}{\mu+c} \right] \nonumber \\
&&+  \left[ \frac{1}{ \mu+\frac{3\lambda-\sqrt{\lambda^2+2^4c^2}}{2^2}}+\frac{1}{ \mu+\frac{3\lambda+\sqrt{\lambda^2+2^4c^2}}{2^2}} \right] \frac{\lambda}{2^2(2 \mu + \frac{3}{2}\lambda)}.
\end{eqnarray}

It is now a simple matter to invert the Laplace transform (\ref{18:55}) and we arrive at the mean hyperbolic distance of the center of mass in an integral form  
\begin{eqnarray} \label{controI}
\lefteqn{E\{\cosh \eta_{cm}(t)\}} \nonumber \\
&=& \left[\frac{1}{2}+\frac{\lambda}{\sqrt{\lambda^2+2^4c^2}}  \right] e^{-t \frac{3\lambda-\sqrt{\lambda^2+2^4c^2}}{2^2}}+\left[\frac{1}{2}-\frac{\lambda}{\sqrt{\lambda^2+2^4c^2}}  \right] e^{-t \frac{3\lambda+\sqrt{\lambda^2+2^4c^2}}{2^2}} \nonumber \\
&&+\frac{\lambda(\lambda+2c)}{2\sqrt{\lambda^2+2^4c^2}} \left[ \int_{0}^{t} e^{cs}e^{-(t-s)\frac{3\lambda-\sqrt{\lambda^2+2^4c^2}}{2^2}} \mathrm{d}s -\int_{0}^{t} e^{cs}e^{-(t-s)\frac{3\lambda+\sqrt{\lambda^2+2^4c^2}}{2^2}} \mathrm{d}s \right]  \nonumber \\
&&+\frac{\lambda(\lambda-2c)}{2\sqrt{\lambda^2+2^4c^2}} \left[ \int_{0}^{t} e^{-cs}e^{-(t-s)\frac{3\lambda-\sqrt{\lambda^2+2^4c^2}}{2^2}} \mathrm{d}s -\int_{0}^{t} e^{-cs}e^{-(t-s)\frac{3\lambda+\sqrt{\lambda^2+2^4c^2}}{2^2}} \mathrm{d}s \right] \nonumber \\
&&+\frac{\lambda}{2^3} \left[ \int_{0}^{t} e^{-\frac{3}{2^2} \lambda s } e^{-(t-s)\frac{3\lambda-\sqrt{\lambda^2+2^4c^2}}{2^2} } \mathrm{d}s +  \int_{0}^{t} e^{-\frac{3}{2^2}\lambda s } e^{-(t-s)\frac{3\lambda+ \sqrt{\lambda^2+2^4c^2}}{2^2} } \mathrm{d}s \right]. 
\end{eqnarray}
The expression (\ref{controI}) can be further developed and simplified by observing that, after some simple calculations, we have that 
 \begin{eqnarray}\label{c2} 
 \lefteqn{\frac{\lambda(\lambda+2c)}{2\sqrt{\lambda^2+2^4c^2}} \left[ \int_{0}^{t} e^{cs}e^{-(t-s)\frac{3\lambda-\sqrt{\lambda^2+2^4c^2}}{2^2}} \mathrm{d}s -\int_{0}^{t} e^{cs}e^{-(t-s)\frac{3\lambda+\sqrt{\lambda^2+2^4c^2}}{2^2}} \mathrm{d}s \right]}  \nonumber \\
&=&\frac{\lambda(\lambda+2c)}{2\sqrt{\lambda^2+2^4c^2}} \frac{2^2}{2^2c+3\lambda-\sqrt{\lambda^2+2^4c^2}} \left[  e^{ct} - e^{-t \frac{3 \lambda-\sqrt{\lambda^2+2^4c^2}}{2^2}}\right] \nonumber \\
&&-\frac{\lambda(\lambda+2c)}{2\sqrt{\lambda^2+2^4c^2}} \frac{2^2}{2^2c+3\lambda+\sqrt{\lambda^2+2^4c^2}} \left[  e^{ct} - e^{-t \frac{3 \lambda+\sqrt{\lambda^2+2^4c^2}}{2^2}}\right].
\end{eqnarray}
Similar manipulations yield
 \begin{eqnarray}\label{c3}
 \lefteqn{\frac{\lambda(\lambda-2c)}{2\sqrt{\lambda^2+2^4c^2}} \left[ \int_{0}^{t} e^{-cs}e^{-(t-s)\frac{3\lambda-\sqrt{\lambda^2+2^4c^2}}{2^2}} \mathrm{d}s -\int_{0}^{t} e^{-cs}e^{-(t-s)\frac{3\lambda+\sqrt{\lambda^2+2^4c^2}}{2^2}} \mathrm{d}s \right]}  \nonumber  \\
&=&\frac{\lambda(\lambda-2c)}{2\sqrt{\lambda^2+2^4c^2}} \;\frac{2^2}{-2^2c+3\lambda-\sqrt{\lambda^2+2^4c^2}} \left[  e^{-ct} - e^{-t \frac{3 \lambda-\sqrt{\lambda^2+2^4c^2}}{2^2}}\right] \nonumber \\
&&-\frac{\lambda(\lambda-2c)}{2\sqrt{\lambda^2+2^4c^2}} \;\frac{2^2}{-2^2c+3\lambda+\sqrt{\lambda^2+2^4c^2}} \left[  e^{-ct} - e^{-t \frac{3 \lambda+\sqrt{\lambda^2+2^4c^2}}{2^2}}\right],
\end{eqnarray} 
and also 
\begin{eqnarray} \label{c1}
\lefteqn{\frac{\lambda}{2^3} \left[ \int_{0}^{t} e^{- \frac{3}{2^2} \lambda s } e^{-(t-s)\frac{3\lambda-\sqrt{\lambda^2+2^4c^2}}{2^2} } \mathrm{d}s +  \int_{0}^{t} e^{- \frac{3}{2^2} \lambda s} e^{-(t-s)\frac{3\lambda+ \sqrt{\lambda^2+2^4c^2}}{2^2} } \mathrm{d}s \right]} \nonumber \\
&=&\frac{\lambda}{2 \sqrt{\lambda^2+2^4c^2}} e^{- \frac{3}{2^2} \lambda t } \left[ e^{t \frac{\sqrt{\lambda^2+2^4c^2}}{2^2}}-e^{-t \frac{\sqrt{\lambda^2+2^4c^2}}{2^2}} \right].
\end{eqnarray}
By inserting results (\ref{c1}), (\ref{c2}), and (\ref{c3}) into (\ref{controI}) we now obtain the final formula 
\begin{eqnarray}
\lefteqn{E\{\cosh \eta_{cm}(t)\}} \nonumber \\
&=&\frac{e^{-\frac{3}{2^2}\lambda t }}{2} \left[  \left(1+ \frac{3\lambda}{\sqrt{\lambda^2+2^4c^2}} \right) e^{\frac{t}{2^2} \sqrt{\lambda^2+2^4c^2}}+\left(1- \frac{3 \lambda}{\sqrt{\lambda^2+2^4c^2}}  \right) e^{-\frac{t}{2^2} \sqrt{\lambda^2+2^4c^2}} \right] \nonumber \\
&&+\frac{\lambda+2c}{2(\lambda+3c)} e^{ct} + \frac{2 \lambda }{\sqrt{\lambda^2+2^4c^2}} e^{-t \frac{3 \lambda-\sqrt{\lambda^2+2^4c^2}}{2^2}} \left[ \frac{\lambda-2c}{2^2c-3\lambda+\sqrt{\lambda^2+2^4 c^2}} -\frac{\lambda+2c}{2^2c+3\lambda-\sqrt{\lambda^2+2^4 c^2}}  \right] \nonumber \\
&&+\frac{\lambda-2c}{2(\lambda-3c)} e^{-ct} + \frac{2 \lambda }{\sqrt{\lambda^2+2^4c^2}} e^{-t \frac{3 \lambda+\sqrt{\lambda^2+2^4c^2}}{2^2}} \left[ \frac{\lambda+2c}{2^2c+3\lambda+\sqrt{\lambda^2+2^4c^2}} -\frac{\lambda-2c}{2^2c-3\lambda-\sqrt{\lambda^2+2^4c^2}}  \right] \nonumber \\
&=& e^{-\frac{3}{2^2}\lambda t } e^{ \frac{t}{2^2}\sqrt{\lambda^2+2^4c^2}} \left[ \frac{1}{2}+\frac{3 \lambda}{2\sqrt{\lambda^2+2^4c^2}}-\frac{2}{\sqrt{\lambda^2+2^4c^2}}\left( \frac{3 \lambda^2- \lambda \sqrt{\lambda^2+2^4c^2} -2^3c^2}{5\lambda-3\sqrt{\lambda^2+2^4c^2}}  \right) \right] \nonumber \\
&&+e^{-\frac{3}{2^2}\lambda t } e^{- \frac{t}{2^2}\sqrt{\lambda^2+2^4c^2}} \left[ \frac{1}{2}-\frac{3 \lambda}{2\sqrt{\lambda^2+2^4c^2}}+\frac{2}{\sqrt{\lambda^2+2^4c^2}}\left( \frac{3 \lambda^2+ \lambda \sqrt{\lambda^2+2^4c^2} -2^3c^2}{5\lambda+3\sqrt{\lambda^2+2^4c^2}}  \right) \right] \nonumber \\
&&+ \frac{\lambda+2c}{2(\lambda+3c)} e^{ct}+\frac{\lambda-2c}{2(\lambda-3c)} e^{-ct} \nonumber \\
&=&\frac{2^3 c^2 e^{-\frac{3}{2^2}\lambda t } }{\sqrt{\lambda^2+2^4c^2}}    \left[ \frac{e^{- \frac{t}{2^2}\sqrt{\lambda^2+2^4c^2}}}{5 \lambda+3 \sqrt{\lambda^2+2^4c^2}} - \frac{e^{ \frac{t}{2^2}\sqrt{\lambda^2+2^4c^2}}}{5 \lambda-3 \sqrt{\lambda^2+2^4c^2}} \right] +\frac{\lambda+2c}{2(\lambda+3c)} e^{ct} +\frac{\lambda-2c}{2(\lambda-3c)} e^{-ct}. \nonumber 
\end{eqnarray} \Fine

\begin{oss}
Apparently a critical point in formula (\ref{12:36}) is $\lambda=3c$. We show that the mean hyperbolic distance (\ref{12:36}) is finite for $\lambda=3c$. We first write $\varepsilon=\lambda-3c$ and $r(\varepsilon)=\sqrt{\varepsilon^2+6 c \varepsilon +5^2 c^2}$ and evaluate the limit 
\begin{eqnarray} \label{13:37}
\lefteqn{\lim_{\lambda \to 3c} \frac{\lambda-2c}{2(\lambda-3c)}e^{-ct}-\frac{2^3 c^2 e^{-\frac{3}{2^2}\lambda t}}{\sqrt{\lambda^2+2^4 c^2}} \frac{e^{\frac{t}{2^2}\sqrt{\lambda^2+2^4c^2}}}{5\lambda-3\sqrt{\lambda^2+2^4c^2}}} \nonumber \\
&=&\lim_{\varepsilon \to 0}  \frac{\varepsilon + c}{2 \varepsilon} e^{-ct}-\frac{2^3 c^2}{r(15c+5 \varepsilon-3r)} e^{-\frac{t}{2^2}(9c+3 \varepsilon-r)},
\end{eqnarray}
which refers to the components of (\ref{12:36}) with diverging coefficients. We expand the second exponential in (\ref{13:37}) as  
\begin{eqnarray}
e^{-\frac{t}{2^2}(9c+3 \varepsilon-r)}&=&e^{-ct}  e^{-\frac{t}{2^2}(5c+3 \varepsilon-r)} =e^{-ct} \left[ 1- \frac{t}{2^2}(5c+3 \varepsilon-r)+\frac{t^2}{2^5}(5c+3 \varepsilon-r)^2+o(\varepsilon^3)  \right], \nonumber
\end{eqnarray}
where we have taken into account that $r(\varepsilon)\to 5c$ for $\varepsilon \to 0$. It is now convenient to write (\ref{13:37}) in the form 
\begin{eqnarray}
e^{-ct} \lim_{\varepsilon \to 0} \left[ \frac{\varepsilon + c}{2 \varepsilon} -\frac{2^3 c^2}{r(15c+5\varepsilon-3r)}  + \frac{2 t c^2 (5c+3 \varepsilon -r)}{r(15 c +5 \varepsilon -3 r)} -\frac{t^2 c^2 (5c+3 \varepsilon -r)^2}{2^2 r(15 c +5 \varepsilon -3 r)} +o(\varepsilon^2) \right]. \nonumber
\end{eqnarray}
By considering that $r(\varepsilon) \sim 5c+\frac{3}{5} \varepsilon$ for $\varepsilon \to 0$, we can easily realize that 
\begin{equation}
\lim_{\varepsilon \to 0} \left[   \frac{\varepsilon + c}{2 \varepsilon} -\frac{2^3 c^2}{r(15c+5\varepsilon-3r)}  \right]=\frac{2\cdot 7}{5^2},\;\;\; \lim_{\varepsilon \to 0} \frac{2 t c^2 (5c+3 \varepsilon -r)}{r(15 c +5 \varepsilon -3 r)}=\frac{3}{2\cdot 5}ct,\;\;\; \lim_{\varepsilon \to 0} \frac{t^2 c^2 (5c+3 \varepsilon -r)^2}{2^2 r(15 c +5 \varepsilon -3 r)}=0. \nonumber
\end{equation}
This suffices to show that near $\lambda=3 c$ the mean hyperbolic distance is finite. Considering also the two additional terms of  (\ref{12:36}) leads to the following asymptotic estimate of the mean hyperbolic distance of the center of mass for large values of $t$ and near $\lambda \sim 3c$. In particular, we have that 
\begin{equation}
E\{\cosh \eta_{cm}(t)\}\sim \frac{5}{12}e^{ct}. \nonumber
\end{equation} 
\begin{figure} 
\centering 
   {\includegraphics[width=15cm, height=6cm]{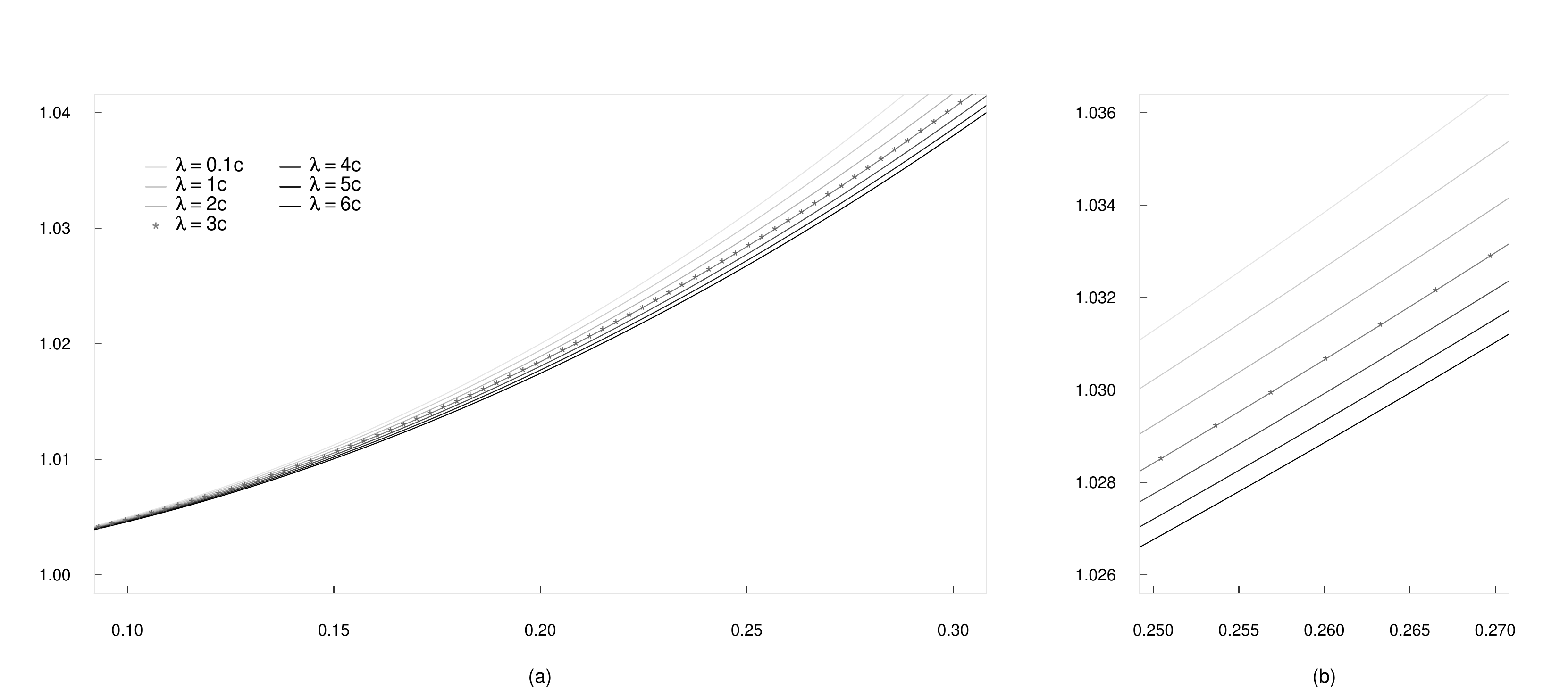}}
 \caption{ In (a) and (b) the mean-value $E\{\cosh \eta_{cm}(t)\}$ is plotted for $c=1$ and different values of $\lambda$, the curve corresponding to $\lambda=3c$ is obtained by plotting the limit value $E\{\cosh \eta_{cm}(t)\}=e^{-ct} \left(  \frac{53}{2^2 5^2}+\frac{3}{10}ct \right)+e^{ct} \frac{5}{2^2 3}+e^{-\frac{7ct}{2} } \frac{2^2 }{5^2 3 }$.} \label{fig1}
 \end{figure}
In Figure \ref{fig1}, the function $E\{\cosh \eta_{cm}(t)\}$ is plotted for $c=1$ and different values of $\lambda$ in the neighborhood of $\lambda=3c$ (including the limiting case $\lambda=3c$).
\end{oss}

\begin{oss} \label{4.2}
We now show that the center of mass (as well as each individual particle) goes further and further away from the origin $O$ of $H_2^+$ (or equivalently it migrates towards the frontier of the Poincar\'e disc $D$) as time $t$ increases. 
In order to show this result we work on (\ref{controI}) and, after some manipulations, we obtain that
\begin{eqnarray} \label{diff}
\frac{\mathrm{d}}{\mathrm{d}t} E\{ \cosh \eta_{cm}(t)\} &=&\frac{2^2 c^2 e^{-\frac{3 }{2^2}\lambda t}}{\sqrt{\lambda^2+2^4c^2}} \sinh \frac{t}{2^2}\sqrt{\lambda^2+2^4c^2} \nonumber \\
&&+  \frac{2 c \lambda^2}{\sqrt{\lambda^2+2^4c^2}} \int_{0}^{t}  e^{-\frac{3 }{2^2} \lambda s} \sinh c(t-s) \sinh \frac{s}{2^2}\sqrt{\lambda^2+2^4c^2}  \mathrm{d}s \nonumber \\
&&+ \frac{2^2 c^2 \lambda}{\sqrt{\lambda^2+2^4c^2}} \int_{0}^{t}  e^{-\frac{3 }{2^2}\lambda s} \cosh c(t-s) \sinh \frac{s}{2^2}\sqrt{\lambda^2+2^4c^2}  \mathrm{d}s.
\end{eqnarray} 
The first term in the right hand side of (\ref{diff}) is produced by the derivatives of all terms of (\ref{controI}) while the integrals stem from the third and fourth term only. From (\ref{diff}) it is easy to check that $\frac{\mathrm{d^2}}{\mathrm{d}t^2} E\{ \cosh \eta_{cm}(t)\}>0$ for all $t$ so that the center of mass gets off from the starting point with positive acceleration.
\end{oss}

\begin{figure}[h] \label{galax}
 \centering
   {\includegraphics[width=15.6cm, height=10.3cm]{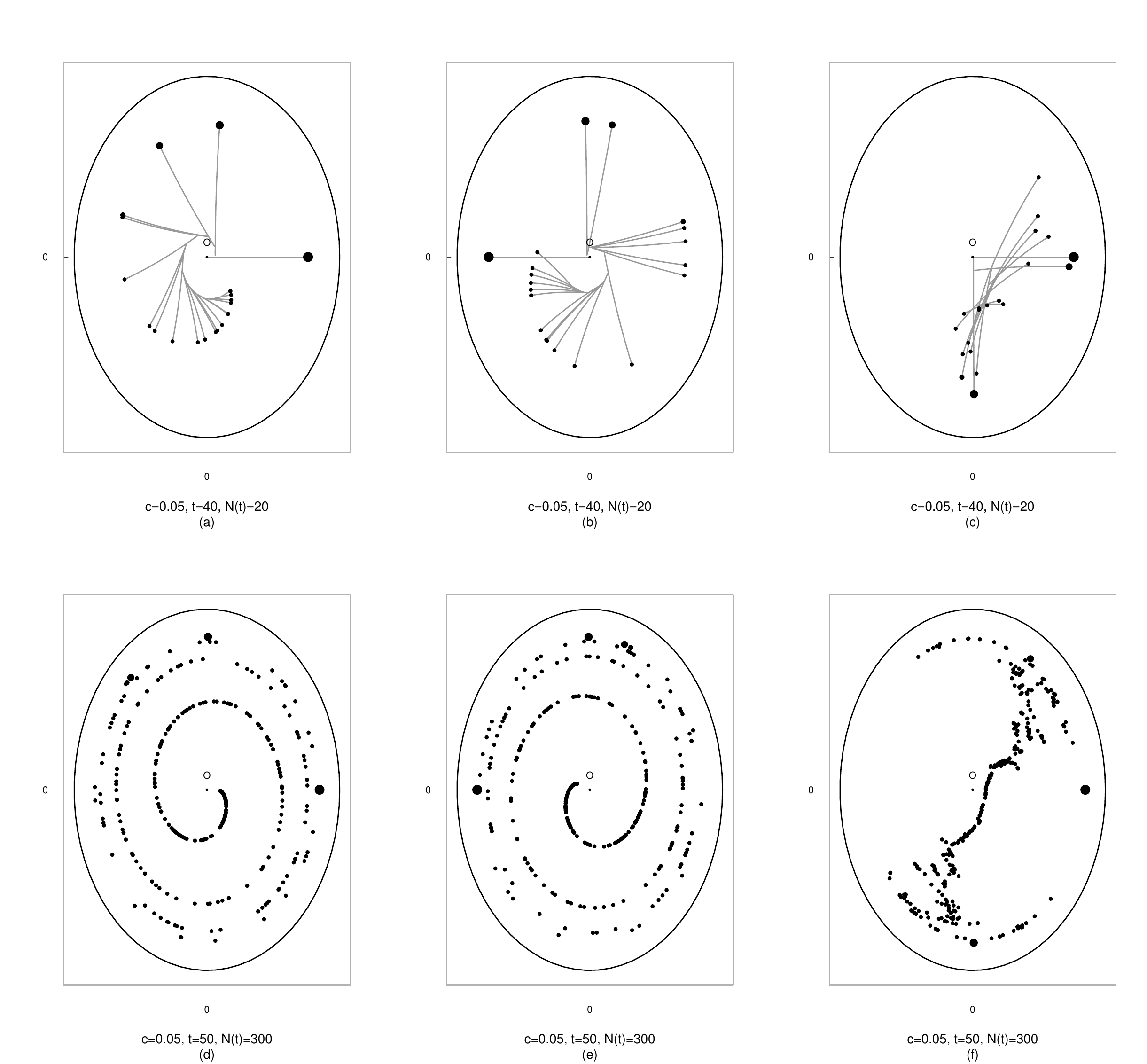}}
 \caption{In (a), (b), and (c) the trajectories and the positions of the splinters at time $t=40$, when $c=0.05$ and $N(t)=20$, are shown. In (d), (e) and (f) only the positions of the splinters at time $t=50$, when $c=0.05$ and $N(t)=300$, are drawn. In (a) and (d) each splinter chooses the clockwise direction, in (b) and (e) each splinter chooses the counterclockwise direction, in (c) and (f) the  clockwise and counterclockwise directions are alternatively chosen. }
 \end{figure}

\section{Equation governing the hyperbolic distance}

We are able to confirm result (\ref{12:36}) by a completely different method based on a system of differential equations governing all the quantities appearing in the mean hyperbolic distance of the center of mass
\begin{equation} \label{14:40}
E\{\cosh \eta_{cm} (t)\}=e^{-\lambda t} \sum_{n=1}^{\infty} \lambda^n \sum_{k=0}^{n-1} \frac{1}{2^{k+1}} G_{n,k}(t)+e^{-\lambda t} \sum_{n=0}^{\infty} \frac{\lambda^n}{2^n} G_{n,n}(t).
\end{equation}   
We start by deriving the difference-differential equations governing the functions $G_{n,k}(t)$, with $t>0$ and $0\le k \le n$. As it has been pointed out in our previous paper \cite{travel}, the $j$-th moment of the conditional hyperbolic distance of a single particle leads to a difference-differential equation of order equal to $j+1$. The same type of phenomenon occurs in the present case where the hyperbolic distance of a cloud of points is envisaged.   

\begin{lem} \label{lemm}
The functions represented by the following multiple integrals 
\begin{equation}
G_{n,k}(t)= \int_{0}^{t} \mathrm{d}s_1 \cdots \int_{s_{k-1}}^{t} \mathrm{d}s_k \; \frac{(t-s_k)}{(n-k)!}^{n-k} \prod_{j=1}^{k} \cosh c(s_j-s_{j-1}) \cosh c(t-s_k), \hspace{1cm}0\le k\le n, \nonumber
\end{equation}
are solutions to the difference-differential equations 
\begin{equation} \label{10:44}
\left\{
\begin{array}{lr} \frac{\mathrm{d}^2 }{\mathrm{d}t^2}G_{n,k}=2 \frac{\mathrm{d}}{\mathrm{d}t}G_{n-1,k} - G_{n-2,k}+ c^2 G_{n,k}, &   k\le n-2,\\
\frac{\mathrm{d}^2 }{\mathrm{d}t^2}G_{n,n-1}=2 \frac{\mathrm{d}}{\mathrm{d}t}G_{n-1,n-1} - G_{n-2,n-2}+ c^2 G_{n,n-1},    & k= n-1,\\
\frac{\mathrm{d}^2 }{\mathrm{d}t^2}G_{n,n}= \frac{\mathrm{d}}{\mathrm{d}t}G_{n-1,n-1} + c^2 G_{n,n},  & k=n.  
\end{array}
\right.
\end{equation}
\end{lem}
\Dim
We start by considering the case $k \le n-2$ where the first derivative reads 
\begin{eqnarray} \label{12:13}
\frac{\mathrm{d}} {\mathrm{d}t}G_{n,k} &=& G_{n-1,k} + c \int_{0}^{t} \mathrm{d}s_1\cdots \int_{s_{k-1}}^{t} \mathrm{d}s_k \frac{(t-s_k)^{n-k}}{(n-k)!} \prod_{j=1}^{k} \cosh c(s_j-s_{j-1}) \sinh c(t-s_k),
\end{eqnarray}
and, from (\ref{12:13}), the second derivative becomes 
\begin{eqnarray}
\frac{\mathrm{d}^2}{\mathrm{d}t^2}G_{n,k}&=&  \frac{\mathrm{d}}{\mathrm{d}t}G_{n-1,k}+ c^2 G_{n,k} \nonumber \\ &&+   c \int_{0}^{t} \mathrm{d}s_1 \cdots \int_{s_{k-1}}^{t} \mathrm{d}s_k \frac{(t-s_k)^{n-k-1}}{(n-k-1)!} \prod_{j=1}^{k} \cosh c(s_j-s_{j-1}) \sinh c(t-s_{k})         \nonumber \\
&=&2\frac{\mathrm{d}}{\mathrm{d}t}G_{n-1,k}-G_{n-2,k} + c^2 G_{n,k}.
\end{eqnarray}
The expression (\ref{12:13}), for $k=n-1$, is qualitatively different and must be handled carefully. In this case we have that   
\begin{eqnarray} \label{14:05}
\frac{\mathrm{d}}{\mathrm{d}t} G_{n,n-1}&=&G_{n-1, n-1}\nonumber \\
&&+ c \int_{0}^{t} \mathrm{d}s_1\cdots \int_{s_{n-2}}^{t} \mathrm{d}s_{n-1} (t-s_{n-1})  \prod_{j=1}^{n-1} \cosh c(s_j-s_{j-1}) \sinh c(t-s_{n-1}),
\end{eqnarray}
and thus, from (\ref{14:05}),  
\begin{eqnarray}
\frac{\mathrm{d}^2}{\mathrm{d}t^2}G_{n,n-1}&=& \frac{\mathrm{d}}{\mathrm{d}t} G_{n-1,n-1} +c^2 G_{n,n-1}\nonumber \\
&&+ c \int_{0}^{t} \mathrm{d}s_1 \cdots \int_{s_{n-2}}^{t} \mathrm{d}s_{n-1} \prod_{j=1}^{n-1} \cosh c (s_j-s_{j-1}) \sinh c (t-s_{n-1}) \nonumber \\
&=& \frac{\mathrm{d}}{\mathrm{d}t}G_{n-1,n-1} +\left( \frac{\mathrm{d}}{\mathrm{d}t}G_{n-1,n-1} - G_{n-2,n-2} \right)+ c^2 G_{n,n-1} \nonumber \\
&=& 2 \frac{\mathrm{d}}{\mathrm{d}t}G_{n-1,n-1} - G_{n-2,n-2}+ c^2 G_{n,n-1}. \nonumber
\end{eqnarray}
We omit the derivation of the last formula in (\ref{10:44}) because it coincides with the result of Lemma 3.1 in \cite{travel}.\Fine

The results of Lemma \ref{lemm} permit us to obtain the equation governing the mean-value of the hyperbolic distance of the center of mass of the randomly moving splinters produced by the disintegration process.    

\begin{teo}
The mean-value of (\ref{dcm}), $E\{\cosh \eta_{cm}(t) \}$, is solution to the non-homogeneus second-order linear equation  
\begin{equation} \label{15:25}
 \frac{\mathrm{d}^2}{\mathrm{d}t^2} u -c^2 u=\frac{\lambda c^2 e^{-\frac{3}{2^2}\lambda t} }{\sqrt{\lambda^2+2^4 c^2}}  \left\{e^{-\frac{t}{2^2} \sqrt{\lambda^2+2^4 c^2}}  -e^{\frac{t}{2^2} \sqrt{\lambda^2+2^4 c^2}} \right\}.
\end{equation}
\end{teo}
\Dim 
We begin by successively deriving the expression (\ref{14:40}) as follows 
\begin{equation}
\frac{\mathrm{d}}{\mathrm{d}t} u=-\lambda u +e^{-\lambda t} \sum_{n=1}^{\infty} \lambda^n \sum_{k=0}^{n-1} \frac{1}{2^{k+1}} \frac{\mathrm{d}}{\mathrm{d}t} G_{n,k} + e^{-\lambda t} \sum_{n=0}^{\infty} \frac{\lambda^n}{2^n} \frac{\mathrm{d}}{\mathrm{d}t} G_{n,n}, \nonumber
\end{equation}
and
\begin{eqnarray} \label{10:53}
\frac{\mathrm{d}^2}{\mathrm{d}t^2}u&=& - 2\lambda \frac{\mathrm{d}}{\mathrm{d}t} u-\lambda^2 u+  e^{-\lambda t} \sum_{n=1}^{\infty}   \lambda^n \sum_{k=0}^{n-1} \frac{1}{2^{k+1}} \frac{\mathrm{d}^2}{\mathrm{d}t^2}G_{n,k} +e^{-\lambda t}  \sum_{n=0}^{\infty}     \frac{\lambda^n}{2^n}  \frac{\mathrm{d}^2}{\mathrm{d}t^2}G_{n,n}.
\end{eqnarray}
We must now insert all the expressions of Lemma \ref{lemm} into (\ref{10:53}) by accurately taking into account the constrains on $n$ and $k$. It is convenient to write the last two terms of (\ref{10:53}) as follows  
\begin{eqnarray} \label{15:04}
 \lefteqn{ \sum_{n=1}^{\infty}   \lambda^n \sum_{k=0}^{n-1} \frac{1}{2^{k+1}} \frac{\mathrm{d}^2}{\mathrm{d}t^2}G_{n,k} + \sum_{n=0}^{\infty}     \frac{\lambda^n}{2^n}  \frac{\mathrm{d}^2}{\mathrm{d}t^2}G_{n,n}} \nonumber \\
 &=&  \sum_{n=2}^{\infty}   \lambda^n \sum_{k=0}^{n-2} \frac{1}{2^{k+1}} \frac{\mathrm{d}^2}{\mathrm{d}t^2}G_{n,k} + \sum_{n=2}^{\infty} \frac{\lambda^n}{2^n}  \frac{\mathrm{d}^2}{\mathrm{d}t^2} G_{n,n-1}+ \frac{\lambda}{2} \frac{\mathrm{d}^2}{\mathrm{d}t^2}G_{1,0} \nonumber \\
&&+ \sum_{n=1}^{\infty}     \frac{\lambda^n}{2^n}  \frac{\mathrm{d}^2}{\mathrm{d}t^2}G_{n,n} +  \frac{\mathrm{d}^2}{\mathrm{d}t^2} G_{0,0} \nonumber \\
 &=&  \sum_{n=2}^{\infty}   \lambda^n \sum_{k=0}^{n-2} \frac{1}{2^{k+1}} \left[ 2 \frac{\mathrm{d}}{\mathrm{d}t}G_{n-1,k} - G_{n-2,k} +c^2 G_{n,k} \right] \nonumber \\  
 &&+ \sum_{n=2}^{\infty} \frac{\lambda^n}{2^n} \left[ 2 \frac{\mathrm{d}}{\mathrm{d}t}G_{n-1, n-1} - G_{n-2,n-2}+c^2 G_{n,n-1} \right] + \lambda \frac{\mathrm{d}}{\mathrm{d}t} G_{0,0} + \frac{\lambda}{2}c^2 G_{1,0} \nonumber \\
 &&+ \sum_{n=1}^{\infty} \frac{\lambda^n}{2^n} \left[  \frac{\mathrm{d}}{\mathrm{d}t}G_{n-1,n-1}+c^2G_{n,n} \right] +  c^2 G_{0,0}. 
 \end{eqnarray}
 By regrouping the above terms we notice that (\ref{15:04}) can take the form 
 \begin{eqnarray} \label{15:17}
\lefteqn{ c^2  \left[ \sum_{n=1}^{\infty} \lambda^n \sum_{k=0}^{n-1} \frac{1}{2^{k+1}}  G_{n,k} +  \sum_{n=0}^{\infty} \frac{\lambda^n}{2^n}  G_{n,n} \right]+ 2 \lambda \frac{\mathrm{d}}{\mathrm{d}t} \left[  \sum_{n=1}^{\infty} \lambda^n \sum_{k=0}^{n-1} \frac{1}{2^{k+1}}  G_{n,k} +  \sum_{n=0}^{\infty} \frac{\lambda^n}{2^n}  G_{n,n} \right]} \nonumber \\
&&-\lambda^2 \left[  \sum_{n=1}^{\infty} \lambda^n \sum_{k=0}^{n-1} \frac{1}{2^{k+1}}  G_{n,k} +  \sum_{n=0}^{\infty} \frac{\lambda^n}{2^n}  G_{n,n} \right] -\frac{\lambda }{2}  \sum_{n=0}^{\infty} \frac{\lambda^{n}}{2^{n}} \frac{\mathrm{d}}{\mathrm{d}t}G_{n,n} +\frac{\lambda^2}{2^2}  \sum_{n=0}^{\infty} \frac{\lambda^{n}}{2^{n}} G_{n,n}\nonumber \\
&=& c^2 e^{\lambda t} u + 2 \lambda \frac{\mathrm{d}}{\mathrm{d}t} [e^{\lambda t} u] - \lambda^2 e^{\lambda t} u  -\frac{\lambda }{2}  \sum_{n=0}^{\infty} \frac{\lambda^{n}}{2^{n}} \frac{\mathrm{d}}{\mathrm{d}t}G_{n,n} +\frac{\lambda^2}{2^2}  \sum_{n=0}^{\infty} \frac{\lambda^{n}}{2^{n}} G_{n,n}.
\end{eqnarray}
Multiplying by $e^{-\lambda t}$ the right hand side of (\ref{15:17}), we have that  
\begin{eqnarray}  \label{15:23}
 \lefteqn{e^{-\lambda t} \sum_{n=1}^{\infty}   \lambda^n \sum_{k=0}^{n-1} \frac{1}{2^{k+1}} \frac{\mathrm{d}^2}{\mathrm{d}t^2}G_{n,k} + e^{-\lambda t} \sum_{n=0}^{\infty}     \frac{\lambda^n}{2^n}  \frac{\mathrm{d}^2}{\mathrm{d}t^2}G_{n,n} }\nonumber \\ 
 &=& c^2 u + 2 \lambda^2 u + 2\lambda  \frac{\mathrm{d}}{\mathrm{d}t} u - \lambda^2 u  -\frac{\lambda }{2}e^{-\lambda t}  \sum_{n=0}^{\infty} \frac{\lambda^{n}}{2^{n}} \frac{\mathrm{d}}{\mathrm{d}t}G_{n,n} +\frac{\lambda^2}{2^2}e^{-\lambda t} \sum_{n=0}^{\infty} \frac{\lambda^{n}}{2^{n}} G_{n,n}.
 \end{eqnarray}
By inserting result (\ref{15:23}) into (\ref{10:53}) we finally obtain 
\begin{eqnarray} \label{kelbert1}
 \frac{\mathrm{d}^2}{\mathrm{d}t^2}u-c^2u&=&- \frac{\lambda }{2} e^{-\lambda t } \sum_{n=0}^{\infty} \frac{\lambda^{n}}{2^{n}} \frac{\mathrm{d}}{\mathrm{d}t}G_{n,n} +\frac{\lambda^2}{2^2}  e^{-\lambda t } \sum_{n=0}^{\infty} \frac{\lambda^{n}}{2^{n}} G_{n,n}.
\end{eqnarray}
The series $e^{-\frac{\lambda t}{2} } \sum_{n=0}^{\infty} \frac{\lambda^{n}}{2^{n}}G_{n,n}$ has been studied in \cite{travel} and represents the mean hyperbolic distance of the particle moving in $H_2^+$ and changing direction at each Poisson event when the rate of the Poisson process is $\lambda/2$. We have that    
\begin{equation}
e^{-\frac{\lambda t}{2} } \sum_{n=0}^{\infty} \frac{\lambda^{n}}{2^{n}}G_{n,n} = \frac{e^{-\frac{\lambda t}{4}}}{2} \left[ \left(1+ \frac{\lambda}{\sqrt{\lambda^2+2^4 c^2}}  \right) e^{\frac{t}{4}\sqrt{\lambda^2+2^4 c^2}}+\left(1- \frac{\lambda}{\sqrt{\lambda^2+2^4 c^2}}  \right) e^{-\frac{t}{4}\sqrt{\lambda^2+2^4 c^2}} \right]. \nonumber
\end{equation} 
Since
\begin{equation}
e^{-\frac{\lambda t}{2}}  \sum_{n=0}^{\infty} \frac{\lambda^{n}}{2^{n}} \frac{\mathrm{d}}{\mathrm{d}t}G_{n,n} =  \frac{\mathrm{d}}{\mathrm{d}t} \left[ e^{-\frac{\lambda t}{2}}  \sum_{n=0}^{\infty} \frac{\lambda^{n}}{2^{n}}G_{n,n}\right]+ \frac{\lambda}{2} e^{-\frac{\lambda t}{2}}  \sum_{n=0}^{\infty} \frac{\lambda^{n}}{2^{n}} G_{n,n},  \nonumber
\end{equation}
the right-hand side of (\ref{kelbert1}) can be rewritten as 
\begin{eqnarray}
\lefteqn{- \frac{\lambda }{2} e^{-\lambda t } \sum_{n=0}^{\infty} \frac{\lambda^{n}}{2^{n}} \frac{\mathrm{d}}{\mathrm{d}t}G_{n,n} +\frac{\lambda^2}{2^2}  e^{-\lambda t } \sum_{n=0}^{\infty} \frac{\lambda^{n}}{2^{n}} G_{n,n}= -\frac{\lambda}{2} e^{-\frac{\lambda t}{2}}  \frac{\mathrm{d}}{\mathrm{d}t} \left[ e^{-\frac{\lambda t}{2}}  \sum_{n=0}^{\infty} \frac{\lambda^{n}}{2^{n}}G_{n,n}\right]} \nonumber \\
&=& -\frac{\lambda}{2} e^{-\frac{\lambda t}{2}} \left\{- \frac{\lambda}{2^3} e^{-\frac{\lambda t}{2^2}} \left[ \left( 1+\frac{\lambda}{\sqrt{\lambda^2+2^4 c^2}} \right) e^{\frac{t}{2^2} \sqrt{\lambda^2+2^4 c^2}}+ \left( 1-\frac{\lambda}{\sqrt{\lambda^2+2^4 c^2}} \right) e^{-\frac{t}{2^2} \sqrt{\lambda^2+2^4 c^2}}\right] \right. \nonumber \\
&&+\left. \frac{ \sqrt{\lambda^2+2^4 c^2}}{2^3} e^{-\frac{\lambda t}{2^2}} \left[ \left( 1+\frac{\lambda}{\sqrt{\lambda^2+2^4 c^2}} \right) e^{\frac{t}{2^2} \sqrt{\lambda^2+2^4 c^2}}- \left( 1-\frac{\lambda}{\sqrt{\lambda^2+2^4 c^2}} \right) e^{-\frac{t}{2^2} \sqrt{\lambda^2+2^4 c^2}}\right] \nonumber  \right\} \nonumber \\
&=& \frac{\lambda^2}{2^4} e^{-\frac{3\lambda t}{2^2}} \left[ \left( 1+\frac{\lambda}{\sqrt{\lambda^2+2^4 c^2}} \right) e^{\frac{t}{2^2} \sqrt{\lambda^2+2^4 c^2}}+ \left( 1-\frac{\lambda}{\sqrt{\lambda^2+2^4 c^2}} \right) e^{-\frac{t}{2^2} \sqrt{\lambda^2+2^4 c^2}}\right] \nonumber \\
&&- \frac{\lambda \sqrt{\lambda^2+2^4 c^2}}{2^4} e^{-\frac{3\lambda t}{2^2}} \left[ \left( 1+\frac{\lambda}{\sqrt{\lambda^2+2^4 c^2}} \right) e^{\frac{t}{2^2} \sqrt{\lambda^2+2^4 c^2}}- \left( 1-\frac{\lambda}{\sqrt{\lambda^2+2^4 c^2}} \right) e^{-\frac{t}{2^2} \sqrt{\lambda^2+2^4 c^2}}\right] \nonumber \\
&=& - e^{-\frac{3\lambda t}{2^2}} e^{\frac{t}{2^2} \sqrt{\lambda^2+2^4 c^2}} \frac{\lambda c^2}{ \sqrt{\lambda^2+2^4 c^2}} + e^{-\frac{3\lambda t}{2^2}} e^{-\frac{t}{2^2} \sqrt{\lambda^2+2^4 c^2}} \frac{\lambda c^2}{ \sqrt{\lambda^2+2^4 c^2}}.   \nonumber
\end{eqnarray}
This lets equation (\ref{15:25}) appear. 
\Fine

The expression (\ref{12:36}) of the mean-value $E\{\cosh \eta_{cm}(t)\}$ derived in Section 4 satisfies the non-homogeneous second-order linear equation (\ref{15:25})  as shown in the next theorem. This proves that the mean hyperbolic distance (\ref{12:36}) is obtained by two different and independent methods.

\begin{teo} \label{morgan}
The mean hyperbolic distance (\ref{12:36}) is the solution to the following Cauchy problem 
\begin{equation} 
\left\{
\begin{array}{lr}  \frac{\mathrm{d}^2}{\mathrm{d}t^2} u -c^2 u=\frac{\lambda c^2 e^{-\frac{3}{2^2}\lambda t} }{\sqrt{\lambda^2+2^4 c^2}}  \left\{e^{-\frac{t}{2^2} \sqrt{\lambda^2+2^4 c^2}}  -e^{\frac{t}{2^2} \sqrt{\lambda^2+2^4 c^2}} \right\},\\
u(0)=1,\\
\left. \frac{\mathrm{d}}{\mathrm{d}t} u(t)\right|_{t=0}=0. \nonumber
\end{array}
\right.
\end{equation}
\end{teo}
\Dim
In order to perform the calculations it is convenient to write  (\ref{12:36}) as 
\begin{equation}
E\{\cosh \eta_{cm}(t)\}=K A(t) B(t) + \frac{\lambda+2c}{2(\lambda+3c)}e^{ct}+\frac{\lambda-2c}{2(\lambda-3c)}e^{-ct} \nonumber
\end{equation}
where 
\begin{equation}
K=\frac{2^3 c^2}{\sqrt{\lambda^2+2^4c^2}},\hspace{0.5cm} A(t)=e^{-\frac{3}{2^2}\lambda t },\hspace{0.5cm} B(t)=\frac{e^{- \frac{t}{2^2}\sqrt{\lambda^2+2^4c^2}}}{5 \lambda+3 \sqrt{\lambda^2+2^4c^2}} - \frac{e^{ \frac{t}{2^2}\sqrt{\lambda^2+2^4c^2}}}{5 \lambda-3 \sqrt{\lambda^2+2^4c^2}}, \nonumber
\end{equation}
and $B(t)$ has derivatives  
\begin{equation}
 B'(t)=-\frac{\sqrt{\lambda^2+2^4c^2}}{2^2} \left[ \frac{e^{- \frac{t}{2^2}\sqrt{\lambda^2+2^4c^2}}}{5 \lambda+3 \sqrt{\lambda^2+2^4c^2}} + \frac{e^{ \frac{t}{2^2}\sqrt{\lambda^2+2^4c^2}}}{5 \lambda-3 \sqrt{\lambda^2+2^4c^2}} \right], \hspace{0.5cm} B''(t)=\frac{\lambda^2+2^4c^2}{2^4} B(t). \nonumber
\end{equation}
Therefore 
\begin{eqnarray}
\lefteqn{\frac{\mathrm{d^2}}{\mathrm{d}t^2} E\{\cosh \eta_{cm}(t)\}-c^2E\{\cosh \eta_{cm}(t)\} =\frac{ \lambda K A}{2} \left[ \frac{5 \lambda}{2^2} B - 3 B' \right] } \nonumber \\
&=&\frac{  2^2 \lambda c^2 e^{-\frac{3}{2^2}\lambda t } }{\sqrt{\lambda^2+2^4c^2}}   \left\{ \frac{5 \lambda}{2^2} \left[ \frac{e^{- \frac{t}{2^2}\sqrt{\lambda^2+2^4c^2}}}{5 \lambda+3 \sqrt{\lambda^2+2^4c^2}} - \frac{e^{ \frac{t}{2^2}\sqrt{\lambda^2+2^4c^2}}}{5 \lambda-3 \sqrt{\lambda^2+2^4c^2}} \right] \right. \nonumber \\
&&+\left.   3 \; \frac{\sqrt{\lambda^2+2^4c^2}}{2^2} \left( \frac{e^{- \frac{t}{2^2}\sqrt{\lambda^2+2^4c^2}}}{5 \lambda + 3 \sqrt{\lambda^2+2^4c^2}} + \frac{e^{ \frac{t}{2^2}\sqrt{\lambda^2+2^4c^2}}}{5 \lambda-3 \sqrt{\lambda^2+2^4c^2}} \right) \right\} \nonumber \\
&=&\frac{\lambda c^2 e^{-\frac{3}{2^2}\lambda t } }{\sqrt{\lambda^2+2^4 c^2}}  \left[ \frac{e^{-\frac{t}{4}\sqrt{\lambda^2+2^4 c^2} }}{5 \lambda+3\sqrt{\lambda^2+2^4 c^2}} (5 \lambda +3 \sqrt{\lambda^2+2^4 c^2}) +\frac{e^{\frac{t}{4}\sqrt{\lambda^2+2^4 c^2} }}{5 \lambda-3\sqrt{\lambda^2+2^4 c^2}} (-5 \lambda +3 \sqrt{\lambda^2+2^4 c^2})  \right] \nonumber \\
&=&\frac{\lambda c^2 e^{-\frac{3}{2^2}\lambda t }}{\sqrt{\lambda^2+2^4 c^2}}   \left[ e^{-\frac{t}{4}\sqrt{\lambda^2+2^4 c^2}}- e^{\frac{t}{4}\sqrt{\lambda^2+2^4 c^2} } \right]. \nonumber
\end{eqnarray}
We can easily check that $E\{\cosh \eta_{cm}(0)\}=1$ and we have that $\left. \frac{\mathrm{d}}{\mathrm{d}t} E\{\cosh \eta_{cm}(t)\} \right|_{t=0} =0$ since  
\begin{eqnarray} 
\lefteqn {\frac{\mathrm{d}}{\mathrm{d}t} E\{\cosh \eta_{cm}(t)\} =K A \left[  -\frac{3 \lambda}{2^2}      B+ B' \right] +c \frac{\lambda+2c}{2(\lambda+3c)} e^{ct}-c \frac{\lambda-2c}{2(\lambda-3c)} e^{-ct}} \nonumber \\
&=& \frac{2 c^2 e^{-\frac{3}{2^2}\lambda t }}{\sqrt{\lambda^2+2^4c^2}}  \left[ \frac{e^{ \frac{t}{2^2}\sqrt{\lambda^2+2^4c^2}}}{5 \lambda-3 \sqrt{\lambda^2+2^4c^2}}(3 \lambda- \sqrt{\lambda^2+2^4c^2}) - \frac{e^{- \frac{t}{2^2}\sqrt{\lambda^2+2^4c^2}}}{5 \lambda+3 \sqrt{\lambda^2+2^4c^2}} (3 \lambda+\sqrt{\lambda^2+2^4c^2})\right] \nonumber \\
&&+c \frac{\lambda+2c}{2(\lambda+3c)} e^{ct}-c \frac{\lambda-2c}{2(\lambda-3c)} e^{-ct}. \nonumber
\end{eqnarray} \Fine 

\begin{oss}
From Theorem \ref{morgan} an interesting relationship between the hyperbolic mean distance of the center of mass and its acceleration can be extracted
\begin{eqnarray}
c^2 E\{ \cosh \eta_{cm}(t)\}-\frac{\mathrm{d}^2}{\mathrm{d}t^2} E\{ \cosh \eta_{cm}(t)\}=\frac{2 \lambda c^2 e^{-\frac{3}{2^2} \lambda t}}{\sqrt{\lambda^2+2^4c^2}} \sinh \frac{t}{2^2} \sqrt{\lambda^2+2^4c^2}.
\end{eqnarray}
For large values of $t$, if $\sqrt{\lambda^2+2^4c^2}-3\lambda>0$, the mean distance and its second derivative tend to coincide while in the opposite case they tend to diverge. In view of Remark \ref{4.2} this means that there is a different rate of growth for the mean distance and its acceleration only if $\sqrt{\lambda^2+2^4c^2}-3\lambda<0$.
\end{oss}

 \begin{figure}[h] \label{D1}
 \centering
   {\includegraphics[width=15.5cm, height=5.4cm]{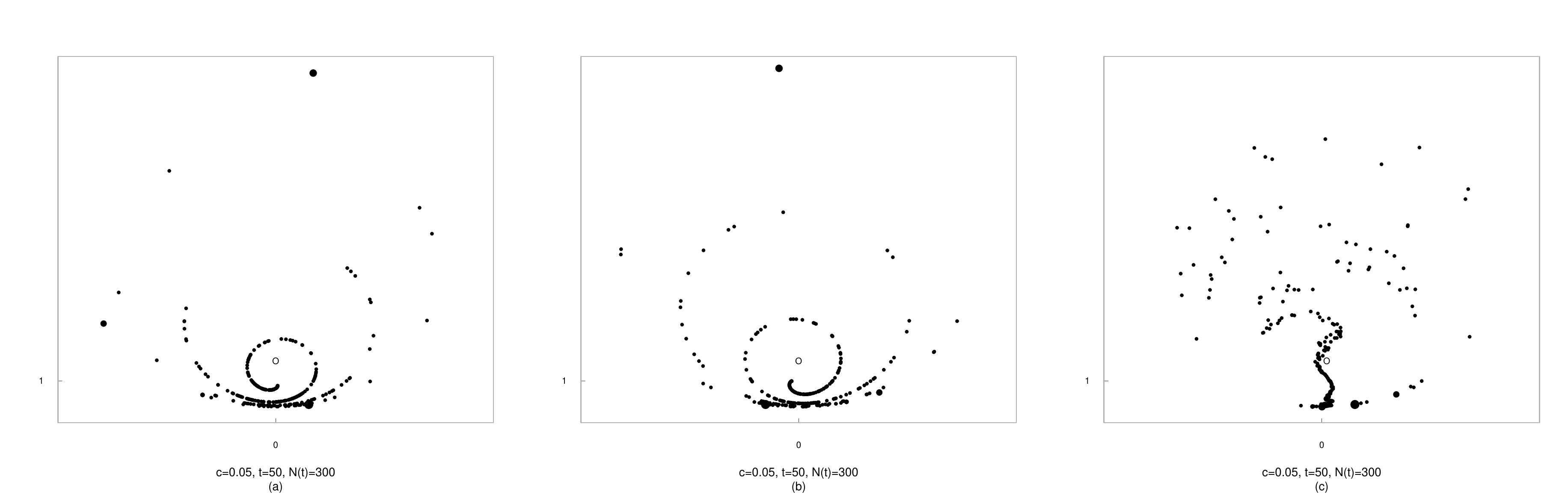}}
 \caption{The position of the splinters in the hyperbolic half-plane $H_2^+$ at time $t=50$, when $c=0.05$ and $N(t)=300$, are drawn. In (a) each splinter chooses the clockwise direction, in (b) each splinter chooses the counterclockwise direction and in (c) the clockwise and counterclockwise directions are alternatively chosen. }
 \end{figure}

\section{Mean hyperbolic distance of the $k$-th splinter}

If $N(t)=n\ge k$, the $k$-th splinter, at time $t$, is located at the hyperbolic distance $\eta_k(t)$ from $O$ equal to 
\begin{equation} \label{14:39}
\cosh \eta_k(t)=\prod_{j=1}^{k+1} \cosh c(S_j-S_{j-1}), \hspace{1cm} 0 \le k \le n,
\end{equation}
The expression (\ref{14:39}) refers to sample paths of particles which change direction until the $k$-th Poisson event and then remain on the same geodesic line until time t. Clearly $\eta_0(t)$ represents the distance of the particle which never changed direction and $\eta_n(t)$ is the distance of the particle that changed direction at all Poisson events. 
  
The mean-value of (\ref{14:39}) becomes 
\begin{eqnarray} \label{14:47}
E\{\cosh \eta_k(t) I_{\{N(t) \ge k\}} \}&=&\sum_{n=k}^{\infty} E\{\cosh \eta_k(t) I_{\{N(t)=n\}}\} \nonumber \\ &=&\sum_{n=k}^{\infty} E\{\cosh \eta_k(t)|N(t)=n\} Pr\{N(t)=n \}.
\end{eqnarray}
We now evaluate the Laplace transform of (\ref{14:47}). By letting $\gamma=\lambda+\mu$ and by considering (\ref{Gk}) we have that 

\begin{eqnarray} \label{18:01}
\int_{0}^{\infty} e^{-\mu t} E\{ \cosh \eta_k(t) I_{\{N(t) \ge k\}} \} \mathrm{d}t &=& \sum_{n=k}^{\infty} \lambda^n \int_{0}^{\infty} e^{-\gamma t} G_{n,k}(t) \mathrm{d}t \nonumber \\
&=&\sum_{n=k}^{\infty} \lambda^n \left( \frac{\gamma}{\gamma^2-c^2} \right)^k \frac{1}{2} \left[ \frac{1}{(\gamma-c)^{n-k+1}}+\frac{1}{(\gamma+c)^{n-k+1}} \right] \nonumber \\
&=&\frac{1}{2}  \left( \frac{\gamma \lambda }{\gamma^2-c^2} \right)^k  \left[ \frac{1}{\gamma-c} \sum_{n=k}^{\infty}  \left( \frac{\lambda}{\gamma-c}\right)^{n-k}+ \frac{1}{\gamma+c} \sum_{n=k}^{\infty}  \left( \frac{\lambda}{\gamma+c}\right)^{n-k}\right] \nonumber \\
&=&\frac{1}{2}  \left( \frac{\gamma \lambda }{\gamma^2-c^2} \right)^k  \left[ \frac{1}{\gamma-c-\lambda}+\frac{1}{\gamma+c-\lambda} \right] \nonumber \\
&=&\frac{1}{2} \left[ \frac{\lambda(\lambda+\mu)}{(\lambda+\mu)^2-c^2} \right]^k  \left( \frac{1}{\mu+c}+\frac{1}{\mu-c} \right) \nonumber \\
&=&\frac{1}{2} \int_{0}^{\infty} e^{-\mu t} \mathrm{d}t \left[ \int_{0}^{t} e^{-c(t-s)} f_k(s) \mathrm{d}s+ \int_{0}^{t} e^{c(t-s)} f_k(s) \mathrm{d}s \right] \nonumber \\
&=& \int_{0}^{\infty} e^{-\mu t} \mathrm{d}t \int_{0}^{t} \cosh c(t-s) f_k(s) \mathrm{d}s 
\end{eqnarray}
where $f_k(s)$ is the inverse Laplace transform of the function $ \left[ \frac{\lambda(\lambda+\mu)}{(\lambda+\mu)^2-c^2} \right]^k$. Let us write  
\begin{eqnarray}
 \left[ \frac{\lambda(\lambda+\mu)}{(\lambda+\mu)^2-c^2} \right]^k&=& \frac{\lambda^k}{2^k} \left[ \frac{1}{\lambda +\mu-c} + \frac{1}{\lambda+\mu+c} \right]^k= \frac{\lambda^k}{2^k} \sum_{r=0}^{k} \binom{k}{r} \frac{1}{(\lambda+\mu-c)^r} \frac{1}{(\lambda+\mu+c)^{k-r}}. \nonumber
\end{eqnarray}
If $k=0$ we have $f_0(s)=\delta(s)$. When $k=1,2,\dots$, recalling that ${\mathit L}\left[ \frac{u^j}{j!}e^{\alpha u} 1_{+}(u) \right]=\frac{1}{(s-\alpha)^{j+1}}$ with $s>\alpha$ and $j=0,1,2 \dots$, we obtain the inverse Laplace transform 
\begin{eqnarray} 
f_k(s)&=&\frac{\lambda^k}{2^{k}} \frac{s^{k-1}}{(k-1)!} e^{s(c - \lambda )}   \nonumber \\
&&+\frac{\lambda^k}{2^k}  \sum_{r=1}^{k-1} \binom{k}{r} \int_{0}^{s} \frac{w^{r-1}}{(r-1)!} e^{w(c-\lambda)} \frac{(s-w)^{k-r-1}}{(k-r-1)!} e^{-(s-w)(\lambda+c)} \mathrm{d}w  \nonumber \\
&&+ \frac{\lambda^k}{2^{k}} \frac{s^{k-1}}{(k-1)!} e^{-s(c + \lambda )}.    \nonumber 
\end{eqnarray}
By inverting the Laplace transform (\ref{18:01}) we obtain for $k=0$ that $E\{\cosh \eta_0(t) I_{\{N(t) \ge 0\}} \}=\cosh ct$, while for $k=1,2, \dots$ we have
\begin{eqnarray} \label{23:53}
\lefteqn{ E\{\cosh \eta_k(t) I_{\{N(t) \ge k\}} \}}  \nonumber \\
&=&\frac{1}{2^k} \int_{0}^{t} \cosh c(t-s) e^{ c s}  \frac{e^{-\lambda s} \lambda^k s^{k-1}}{(k-1)!}   \mathrm{d}s \nonumber \\
&&+\frac{1}{2^k}  \sum_{r=1}^{k-1} \binom{k}{r}  \int_{0}^{t} \cosh c(t-s) e^{-\lambda s} {\lambda^k} s^{k-1}  \int_{0}^{1} e^{cs(2y-1)} \frac{y^{r-1}  (1-y)^{k-r-1}}{(r-1)!  (k-r-1)!}  \mathrm{d}y \;  \mathrm{d}s \nonumber \\
&&+\frac{1}{2^k} \int_{0}^{t} \cosh c(t-s) e^{-c s}  \frac{e^{-\lambda s} \lambda^k s^{k-1}}{(k-1)!}   \mathrm{d}s \nonumber \\
&=&\frac{1}{2^k} \int_{0}^{t} \cosh c(t-s) e^{ c s}  \frac{e^{-\lambda s} \lambda^k s^{k-1}}{(k-1)!}   \mathrm{d}s \nonumber \\
&&+\frac{1}{2^k}  \sum_{r=1}^{k-1} \binom{k}{r}  \int_{0}^{t} \cosh c(t-s)  E_{Y_{r,k}} \{ e^{cs(2Y_{r,k}-1)}\}  \frac{e^{-\lambda s} {\lambda^k} s^{k-1}}{\Gamma(k)}   \mathrm{d}s \nonumber \\
&&+\frac{1}{2^k} \int_{0}^{t} \cosh c(t-s) e^{-c s}  \frac{e^{-\lambda s} \lambda^k s^{k-1}}{(k-1)!}   \mathrm{d}s \nonumber \\
&=&\frac{1}{2^k} \int_{0}^{t} \cosh c(t-s) \left[  e^{ c s} + \sum_{r=1}^{k-1} \binom{k}{r}  E_{Y_{r,k}} \{e^{cs(2Y_{r,k}-1)} \} +e^{-c s}  \right] g(s; k,\lambda) \; \mathrm{d}s 
 \end{eqnarray}
 where $g(s; k,\lambda) =  \frac{e^{-\lambda s} {\lambda^k} s^{k-1}}{\Gamma(k)} $, $Y_{r,k} \sim \mathrm{Beta}(r,k-r)$ and $(2Y_{r,k}-1)\in(-1,1)$. If we adopt the convention that $Y_{0,k}=1$ and $Y_{k,k}=-1$ it is possible to simplify the last expression rewriting it as    
\begin{eqnarray}
E\{\cosh \eta_k(t) I_{\{N(t) \ge k\}} \}&=&\frac{1}{2^k} \int_{0}^{t} \cosh c(t-s)  \sum_{r=0}^{k} \binom{k}{r}  E_{Y_{r,k}} \{e^{cs(2Y_{r,k}-1)} \}  g(s; k,\lambda) \; \mathrm{d}s \nonumber \\
&=& \frac{1}{2^k} \int_{0}^{t} \cosh c(t-s) \; h(k,c,s) \; g(s; k,\lambda) \; \mathrm{d}s \nonumber
\end{eqnarray} 
 where $h(k,c,s)= \sum_{r=0}^{k} \binom{k}{r}  E_{Y_{r,k}} \{e^{cs(2Y_{r,k}-1)} \}$. The last result suggests regarding the mean hyperbolic distance of a particle generated at the $k$-th Poisson event as the mean hyperbolic distance of a particle moving on the first geodesic line and stopping at a randomly distributed time; the law of the stopping time is a suitable combination of a Beta and Gamma distributions.

\begin{oss} In some particular cases it is possible to check easily that formula (\ref{23:53}) gives exactly the same expression of the mean hyperbolic distance of the $k$-th particle as that obtained starting from formula (\ref{14:47}). In particular for $k=1$ we have   
\begin{eqnarray}
\lefteqn{E\{\cosh \eta_1(t) I_{\{N(t) \ge 1\}} \}} \nonumber \\
&=&\frac{1}{2} \int_{0}^{t} \cosh c(t-s) \left[  e^{ c s}   +e^{-c s}  \right] g(s; 1,\lambda) \mathrm{ds}  = \int_{0}^{t} \cosh c(t-s) \cosh cs\; \lambda e^{-\lambda s} \mathrm{d}s. \nonumber
\end{eqnarray}
For $k=2$ we obtain 
\begin{eqnarray}
\lefteqn{E\{\cosh \eta_2(t) I_{\{N(t) \ge 2\}} \}} \nonumber \\
&=& \frac{1}{2^2} \int_{0}^{t} \cosh c(t-s) \left[  e^{ c s} + 2 E_{Y_{1,2} }\{ e^{cs(2Y_{1,2} -1)}\}  +e^{-c s}  \right] g(s; 2,\lambda) \; \mathrm{d}s \nonumber \\
&=&  \frac{1}{2} \int_{0}^{t}   \cosh c(t-s)  \cosh cs  \frac{\lambda^2 s e^{-\lambda s}}{\Gamma(2)} \mathrm{d}s + \frac{1}{2 }  \int_{0}^{t} \cosh c(t-s) \; E_{Y_{1,2} }\{ e^{cs(2Y_{1,2} -1)}\} \; \frac{\lambda^2 s e^{-\lambda s}}{\Gamma(2)}\; \mathrm{d}s \nonumber 
\end{eqnarray}
where $Y_{1,2} \sim \mathrm{Beta}(1,1)$.
\end{oss}

\end{document}